%% file: nondivergence.tex
\documentclass[10pt]{article}

\usepackage[latin1]{inputenc}
\usepackage[T1]{fontenc}
\usepackage[francais]{babel}
\usepackage{amsfonts,amssymb,graphicx,psfrag}
\usepackage{euscript,epsfig}

\newtheorem{theo}{Th\'eor\`eme}[section]

\newtheorem{prop}[theo]{Proposition}
\newtheorem{corol}[theo]{Corollaire}
\newtheorem{lem}[theo]{Lemme}

\newenvironment{demon}[1]{{\flushleft \bf D\'emonstration #1: }}{\hfill $\square$ \vspace{5mm}}
\newenvironment{demo}{{\flushleft \bf D\'emonstration: }}{\hfill $\square$ \vspace{5mm}}

\newcommand{\R}{\mathbf R}
\renewcommand{\H}{\mathbb H}
\newcommand{\N}{\mathbf N}
\newcommand{\M}{\widetilde{M}}
\newcommand{\D}{\partial\widetilde{M}}
\newcommand{\DD}{{\partial}^2\widetilde{M}}
\newcommand{\DDR}{{\partial}^2\widetilde{M}\times\R}
\newcommand{\Z}{\mathbf Z}

\newcommand{\tm}{T^1M}
\newcommand{\ttm}{T^1\widetilde{M}}

\begin{document}

\title{ Lemme de l'ombre et non divergence des horosph\`eres d'une vari\'et\'e g\'eom\'etriquement finie}

\author{ Barbara SCHAPIRA\\
\\
\em  MAPMO,
Universit\'e d'Orl\'eans,
Rue de Chartres,\\
\em BP 6759,
45067 Orl\'eans cedex 2, France\\
\em schapira@labomath.univ-orleans.fr\\}

\maketitle

%

\begin{abstract}
Dans ce travail, nous d\'emontrons le Lemme de l'Ombre sur des vari\'et\'es g\'eom\'etriquement finies \`a courbure n\'egative variable.
Nous en d\'eduisons un r\'esultat de non divergence des horosph\`eres de telles vari\'et\'es.
\end{abstract}


\section{Introduction}

Soit $M$ une vari\'et\'e riemannienne \`a courbure n\'egative. 
Le flot g\'eod\'esique agissant sur son fibr\'e unitaire tangent $\tm$ est 
un flot hyperbolique; en particulier, $\tm$ est feuillet\'e par les
 vari\'et\'es fortement instables de ce flot, encore appel\'ees {\it horosph\`eres fortement instables}.

Si $M$ est une surface hyperbolique, les feuilles de ce feuilletage 
sont les orbites du flot horocyclique $(h^t)_{t\in\R}$ agissant sur $\tm$. 
On sait que lorsque $M$ est compacte, ce flot est minimal: les orbites sont toutes denses dans $\tm$.
Quand elle est seulement de volume fini, certaines orbites sont p\'eriodiques, 
les autres sont denses dans $\tm$ (Hedlund \cite{Hedlund}). 
En particulier, leur projection sur $M$ revient infiniment souvent dans la partie compacte de la vari\'et\'e.
Ce r\'esultat qualitatif de retour dans un compact a \'et\'e \'etendu par Margulis \cite{Margulis} \`a 
l'action de flots unipotents sur l'espace des r\'eseaux $\Lambda_n=SL(n,\R)/SL(n,\Z)$.

Ce fut  Dani qui, le premier, pr\'ecisa ce comportement en une estim\'ee quantitative: 
dans le cas du flot horocyclique d'une surface hyperbolique dans \cite{dani2bis}, puis dans le cas d'un flot unipotent sur $\Lambda_n$ dans \cite{dani3}, 
il montre que pour tout  $\epsilon>0$, il existe un compact $K_{\epsilon}$ de $\tm$ tel 
qu'une orbite non p\'eriodique de $(h^t)_{t\in\R}$ de 
longueur $T$ suffisamment grande passe un temps sup\'erieur \`a 
 $(1-\epsilon)T$ dans $K_{\epsilon}$.
Rappelons que ce r\'esultat, dit de non divergence, est un outil 
essentiel dans le th\'eor\`eme d'\'equidistribution des flots unipotents de Ratner.

Dans un cadre diff\'erent, mentionnons encore le r\'esultat de Minsky et Weiss \cite{MW} de non 
divergence des orbites du flot horocyclique de Teichm\"uller.

Le but du pr\'esent travail est d'\'etablir l'analogue de ce ph\'enom\`ene 
de non divergence en l'\'etendant dans plusieurs directions: 
nous consid\'erons des vari\'et\'es qui sont d'une part de dimension quelconque, 
d'autre part de courbure n\'egative variable, et enfin qui sont g\'eom\'etriquement finies de volume infini.

Nous \'etablissons un r\'esultat de ce type sur des vari\'et\'es 
de volume infini \`a courbure variable major\'ee par $-1$.
Dans ce cadre plus g\'en\'eral, nous restreignons notre \'etude \`a 
l'ensemble non errant ${\cal E}$ du feuilletage horosph\'erique de $\tm$.
On sait alors (Dal'bo \cite{Dalbo}) que certaines feuilles de ${\cal E}$  sont compactes, 
et les autres sont denses dans ${\cal E}$. 
 En particulier, leur projection sur $M$ revient infiniment 
souvent dans la partie compacte de la vari\'et\'e. 
Notons ${\cal E}_R\subset{\cal E}$ l'union de ces feuilles denses dans ${\cal E}$.

En l'absence de param\'etrage naturel de ces feuilles par un flot, 
nous consid\'erons sur chaque feuille de grandes boules not\'ees $B^+(u,r)$, avec $r>0$, 
pour une distance adapt\'ee.

Pour donner un sens au << temps pass\'e >> par une feuille dans un compact, 
nous \'etudions le comportement de moyennes sur ces boules 
pour une mesure, not\'ee $\mu_{H^+}$ sur chaque feuille $H^+$, construite \`a partir d'une mesure naturelle 
sur le bord $\D$ du rev\^etement universel $\M$ de $M$, la {\it mesure de Patterson}.

Plus pr\'ecis\'ement, si $u\in\tm$ et $r>0$, nous d\'efinissons 
la moyenne d'une fonction $\psi:\tm\to\R$ continue par:
$$
M_{r,u}(\psi):=\oint_{B^+(u,r)}\psi(v)\,d\mu_{H^+}(v).
$$

Sous une condition de divergence des cusps de $M$, not\'ee $(*)$ et 
\'enonc\'ee plus loin, notre r\'esultat principal est le suivant:\\

\noindent
{\bf Th\'eor\`eme \ref{nondivergence}:} {\em 
Soit $M$ une vari\'et\'e g\'eom\'etriquement finie de 
courbures sectionnelles major\'ees par $-1$ dont les cusps v\'erifient 
la condition $(*)$.
Soit $\varepsilon>0$ fix\'e, et $C\subset \tm$ un compact de l'ensemble non errant du feuilletage horosph\'erique. 
Il existe un compact $K_{\varepsilon,C}\subset{\cal E}$, tel que pour 
tout vecteur non errant $u$ de ${\cal E}_R\cap C$  et pour tout $r\ge 0$, on a:
$$
M_{r,u}(K_{\varepsilon,C})\ge 1-\varepsilon
$$
}

Ce r\'esultat assure en particulier que les moyennes $(M_{r,u})_{r\ge 0}$ 
ne tendent pas faiblement vers $0$ quand $r\to\infty$.
Autrement dit, il n'y a aucune perte de masse dans les cusps.

Notons que  ce r\'esultat de non divergence a \'et\'e d\'emontr\'e par Rudolph \cite{Rudolph} 
en courbure constante, mais seulement  pour presque tout $u\in{\cal E}_r$ (pour la mesure de Patterson-Sullivan), 
ce qui lui a permis d'en d\'eduire l'ergodicit\'e du feuilletage.\\

Pour \'etablir ce r\'esultat, il nous faudra donner une estim\'ee 
pr\'ecise de la mesure de Patterson de petits ouverts du bord, les {\it ombres}.
Ce r\'esultat, appel\'e {\it Lemme de l'Ombre}, est originellement d\^u \`a Sullivan sur les vari\'et\'es hyperboliques compactes, 
et a \'et\'e \'etendu \`a toutes les vari\'et\'es
compactes ou convexes-cocompactes \`a courbure n\'egative variable, permettant alors 
de r\'einterpr\'eter la mesure de Patterson comme la mesure de Hausdorff 
de l'ensemble limite $\Lambda_{\Gamma}$ de $\Gamma$ dans $\D$.
Stratmann et Velani \cite{SV} l'ont g\'en\'eralis\'e sur les vari\'et\'es hyperboliques 
g\'eom\'etriquement finies, et nous adaptons une preuve de leur r\'esultat 
due \`a Peign\'e \cite{peigne} pour l'\'etablir sur les vari\'et\'es 
g\'eom\'etriquement finies de courbure variable inf\'erieure \`a $-1$ 
qui satisfont la condition $(*)$ \'enonc\'ee ci-dessous.

Pr\'ecisons d'abord quelques notations. Plut\^ot que les ombres consid\'er\'ees 
par Sullivan, nous consid\'ererons des ensembles comparables: 
si $o$ est un point de $\M$, $\xi$ un point du bord $\D$, et $t\ge 0$, 
nous noterons $V(o,\xi,t)$ l'ensemble des points du
 bord dont le projet\'e sur le rayon $[o\xi)$ est \`a distance au moins $t$ de $o$. 

Sur une vari\'et\'e g\'eom\'etriquement finie, l'ensemble int\'eressant d'un point 
de vue dynamique se d\'ecompose en une partie compacte et un nombre fini de pointes, 
les {\em cusps}. Pour simplifier, nous supposerons qu'elle n'en a qu'un seul. 
A ce cusp est associ\'ee une classe de conjugaison de sous-groupes 
paraboliques de $\Gamma$, i.e. des sous-groupes de $\Gamma$ fixant exactement un 
point de $\Lambda_{\Gamma}$, et maximaux pour cette propri\'et\'e.

Nous aurons besoin de faire l'hypoth\`ese suivante, not\'ee $(*)$, 
sur la croissance des sous-groupes paraboliques de $\Gamma$. 
Si $o\in\M$ est un point fix\'e, pour tout sous-groupe parabolique $\Pi$ de $\Gamma$, 
il existe une constante $D\ge 1$ telle que 
\begin{eqnarray*}
\quad\quad\quad\frac{1}{D}\exp(\delta_{\Pi}T)\,\le\,
\sharp \{p\in\Pi,d(o, po)\in[T,T+1[\}\,\le\,D\exp(\delta_{\Pi}T),\quad\quad\quad (*)
\end{eqnarray*}
o\`u pour tout sous-groupe $G$ de $\Gamma$, $\delta_G$ d\'esigne 
l'{\it exposant critique} de $G$, d\'efini par:
$$
\delta_G=\limsup_{T\to\infty}\,\frac{1}{T}\log\,\sharp\{g\in G,\, d(o,go)\le T\}
$$
En particulier, $\delta_G\le \delta_{\Gamma}$.
Nous verrons plus loin (proposition \ref{satisfaite}) que cette hypoth\`ese est toujours satisfaite 
si $M$ est une vari\'et\'e localement sym\'etrique de rang $1$.

Une {\it densit\'e conforme invariante} par $\Gamma$ de dimension $\delta>0$ 
est une famille $\nu=(\nu_x)_{x\in\M}$ de 
mesures \'equivalentes sur le bord $\D$, 
telles que pour tout $\gamma\in\Gamma$ et $x\in\M$, $\nu_{\gamma x}=\gamma_*\nu_x$, et pour tous $(x,y)\in\M^2$, 
$d\nu_x(\xi)=\exp(-\delta\beta_{\xi}(x,y))\,d\nu_y(\xi)$, o\`u $\beta_{\xi}(x,y)$ d\'esigne le cocycle de Busemann (voir paragraphe \ref{s21})

Le Lemme de l'Ombre s'\'enonce alors:\\

\noindent
{\bf Th\'eor\`eme \ref{courbvariable}:} {\it 
Soit $M=\M/\Gamma$ une vari\'et\'e g\'eom\'etriquement finie \`a un cusp $C_1$ qui satisfait l'hypoth\`ese $(*)$. 
Alors pour toute  densit\'e conforme $\nu=(\nu_x)_{x\in\M}$ $\Gamma$-invariante sans atomes normalis\'ee 
de dimension $\delta$ et de support 
$\Lambda_{\Gamma}$, il existe des constantes $A_0>0$ et $A_1>0$ 
telles que pour tout $\xi\in\Lambda_{\Gamma}$ et $t\ge 0$, en notant $\xi(t)$ le point \`a distance $t$ de $o$ sur $[o\xi)$:\\
{\bf a- } si $\xi(t)$ appartient \`a un relev\'e de la partie compacte sur $\M$, alors
$$\quad\frac{1}{A_0}\,e^{-\delta t}\le\nu_{o}(V(o,\xi,t))\,\le A_0 \,e^{-\delta t},$$
{\bf b- } si $\xi(t)$ est dans un relev\'e du cusp, alors:
$$
\frac{1}{A_1}\,e^{-\delta t+(2\delta_{\Pi}-\delta)\,d(\xi(t),\Gamma o)}\,\le\,
\nu_{o}(V(o,\xi,t))\,\le \,A_1 \,e^{-\delta t+(2\delta_{\Pi}-\delta)\,d(\xi(t),\Gamma o)} .
$$
}

\noindent 
L'organisation du texte est la suivante: nous commen\c cons 
(section \ref{s2}) par une introduction aux vari\'et\'es 
g\'eom\'etriquement finies, et nous prouvons une s\'erie de 
lemmes g\'eom\'etriques \'el\'ementaires utiles dans la suite.
Au paragraphe \ref{s3}, nous \'enon\c cons et d\'emontrons le 
Lemme de l'Ombre (th\'eor\`eme \ref{courbvariable}). 
Dans la derni\`ere partie (section \ref{s4}), apr\`es des 
rappels sur le feuilletage fortement instable et les moyennes qui 
nous int\'eressent, nous prouvons le th\'eor\`eme \ref{nondivergence}.

Je remercie ma directrice de th\`ese Martine Babillot pour m'avoir fait 
d\'ecouvrir de belles math\'ematiques pendant ma th\`ese, et pour ses lectures, commentaires et corrections 
 de nombreuses versions pr\'eliminaires de ce travail.


\section{Vari\'et\'es g\'eom\'etriquement finies}\label{s2}

\subsection{G\'en\'eralit\'es}\label{s21}

Nous renvoyons \`a Bowditch \cite{bowditch} pour un expos\'e complet 
sur les vari\'et\'es g\'eom\'etriquement finies \`a courbure n\'egative pinc\'ee, 
et \`a Roblin \cite{Roblin} pour des  compl\'ements dans le cas des 
vari\'et\'es \`a courbure seulement major\'ee; 
nous n'exposons ici que ce dont nous avons besoin.

Soit $M$ une vari\'et\'e riemannienne \`a courbure sectionnelle major\'ee par $-1$,
 $\widetilde{M}$ son rev\^etement universel, et $\Gamma=\pi_1(M)$ son groupe fondamental. 
Notons $\tm$ le fibr\'e unitaire tangent de $M$, et $\pi:\tm\to M$ la projection canonique. 
Nous noterons $d$ la distance riemannienne sur $M$ et $\M$.

Le {\it bord \`a l'infini} $\D$ de $\M$ permet de compactifier $\M$ en $\overline{M}=\M\cup\D$. 
Le groupe $\Gamma$ agit sur $\M$ par isom\'etries, et sur $\D$ par hom\'eomorphismes.
L'{\it ensemble limite} $\Lambda_{\Gamma}\subset \D$ de $\Gamma$ est le plus petit ferm\'e $\Gamma$-invariant de $\D$. 
C'est aussi l'ensemble des points d'accumulation dans $\D$ de l'orbite d'un point quelconque $o\in\M$ par $\Gamma$: 
$\Lambda_{\Gamma}=\overline{\Gamma o}\setminus \Gamma o$.

Le {\em  flot g\'eod\'esique }$g=(g^t)_{t\in\R}$ de $M$ agit sur $\tm$ 
 en associant \`a un vecteur $v$ le vecteur $\dot{c}_v(t)$ 
tangent \`a l'unique g\'eod\'esique $(c_v(t))_{t\in\R}$ telle que $\dot{c}_v(0)=v$. Il se rel\`eve sur $\ttm$ en le 
flot g\'eod\'esique de $\M$, not\'e $\tilde{g}=(\tilde{g}^t)_{t\in\R}$.
L'{\em ensemble non errant} $\Omega\subset\tm$ du flot g\'eod\'esique s'identifie (Eberlein, \cite{Eberlein}) 
\`a l'ensemble 
des vecteurs $v\in\tm$ dont un relev\'e $\widetilde{v}\in\ttm$ d\'efinit une g\'eod\'esique dont les 
deux extr\'emit\'es sont dans $\Lambda_{\Gamma}$.

Un point $\xi$ de l'ensemble limite $\Lambda_{\Gamma}$ est dit {\it radial} s'il existe un point $o\in\M$
et une infinit\'e de points de l'orbite $\Gamma o$ \`a distance born\'ee du rayon $[o\xi)$. L'ensemble des
points limite radiaux sera not\'e $\Lambda_R\subset\Lambda_{\Gamma}$.

Si $\xi\in\Lambda_{\Gamma}$ est l'unique point fixe d'une isom\'etrie parabolique de $\Gamma$, 
il est dit {\em parabolique}. Le stabilisateur dans $\Gamma$ d'un tel point sera appel\'e un {\em sous-groupe parabolique maximal}.
Le point $\xi$ est {\em parabolique born\'e} si son stabilisateur $\Pi\subset\Gamma$ agit de mani\`ere
cocompacte sur $\Lambda_{\Gamma}\setminus\{\xi\}$. 
Nous noterons $\Lambda_{pb}$ l'ensemble des points paraboliques born\'es de $\Lambda_{\Gamma}$.

Le groupe $\Gamma$ est dit {\em cocompact} si $M$ est compacte, 
ce qui implique que $\Omega=\tm$ est compact, et $\Lambda_{\Gamma}=\Lambda_R=\D$.
Il est dit {\em convexe-cocompact} si $\Omega$ est compact, et dans ce cas $\Lambda_{\Gamma}=\Lambda_R$.
Enfin, il est {\em g\'eom\'etriquement fini} si $\Lambda_{\Gamma}=\Lambda_R\cup\Lambda_{pb}$.

Le {\em cocycle de Busemann} est d\'efini sur $\D\times\M\times\M$ par:
$$\beta_{\xi}(x,y)=\lim_{z\to\xi}d(x,z)-d(y,z)="d(x,\xi)-d(y,\xi)"$$
C'est une fonction continue qui v\'erifie la relation de cocycle: $\beta_{\xi}(x,y)+\beta_{\xi}(y,z)=\beta_{\xi}(x,z)$.
Ce cocycle permet de donner des coordonn\'ees sur $\ttm$. 
Si $u\in\ttm$, on notera $u^+\in\D$ (resp. $u^-$) 
l'extr\'emit\'e $c_u(+\infty)$ (resp. $c_u(-\infty)$) de l'unique g\'eod\'esique $c_u$ telle que $\dot{c}_u(0)=u$. 
L'ensemble des g\'eod\'esiques orient\'ees de 
$\ttm$ est en bijection avec le << double bord >> $\DD:=\D\times\D\setminus\{(\xi,\xi),\,\xi\in\D\}$. 
Soit $o\in\M$ un point fix\'e une fois pour toutes.
Alors l'application ci-dessous est un hom\'eomorphisme:
\begin{eqnarray*}
   \ttm & \to & \DDR \\
 v & \mapsto & (v^-, v^+,\beta_{v^-}(\pi(v),o))
\end{eqnarray*}

\begin{figure}[!ht]
\begin{center}
\input{coordo_bordND.pstex_t}
\caption{Coordonn\'ees sur $\ttm$}
\end{center} 
\end{figure}
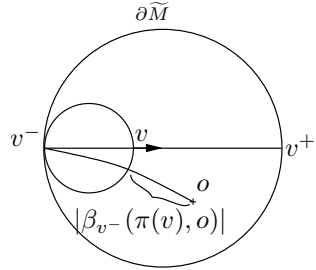

Sur $\DDR$, les actions de $\Gamma$ et de $\tilde{g}$ commutent et s'\'ecrivent:
$$
\gamma\,(u^-,u^+,s)=(\gamma u^-,\gamma u^+,s+\beta_{u^-}(o,\gamma^{-1}o))\quad\mbox{et}\quad 
\tilde{g}^t(u^-,u^+,s)=(u^-,u^+,s+t)
$$
Ainsi on a aussi des hom\'eomorphismes $\tm\simeq(\DDR)/\Gamma$ et $\Omega\simeq(\Lambda_{\Gamma}^2\times\R)/\Gamma$. 

Une {\em horosph\`ere} $H\subset\M$ centr\'ee en $\xi$ est une ligne de niveau de l'application $y\to\beta_{\xi}(y,o)$. 
Une {\em horoboule} ${\cal H}\subset\M$ centr\'ee en $\xi$ est un sous-ensemble 
${\cal H}=\{y\in\M, \beta_{\xi}(y,o)\leq C\}$, avec $C\in\R$.
C'est un ensemble g\'eod\'esiquement convexe de $\M$.

Notons $\widetilde{C}(\Gamma)$ l'enveloppe convexe dans $\M$ de 
l'ensemble limite $\Lambda_{\Gamma}$. 
Le {\em coeur de Nielsen} $N_{\Gamma}$ de $M$ (ou de $\Gamma$) 
est le quotient $N_{\Gamma}=\widetilde{C}(\Gamma)/\Gamma$.
La vari\'et\'e $M$ est g\'eom\'etriquement finie si et seulement si 
il se d\'ecompose en une union finie (voir Bowditch \cite{bowditch} en courbure pinc\'ee, 
et Roblin \cite{Roblin} Prop. 1.10 dans le cas g\'en\'eral):
$$
 N_{\Gamma}\,=\,C_0\sqcup  C_1\dots\sqcup C_k,
$$
o\`u $C_0$ est un ensemble relativement compact, 
de diam\`etre not\'e $\Delta$, et les $C_l$, $1\le l\le k$ (en nombre fini) 
sont les {\em cusps}: pour chaque $l$, $C_l$ est isom\'etrique au 
quotient de ${\cal H}^l\cap \widetilde{C}(\Gamma)$ par $\Pi_l$, o\`u ${\cal H}^l$
est une horoboule centr\'ee en un point parabolique $\xi^l$, 
et $\Pi_l$ est le stabilisateur de $\xi^l$. 
Nous noterons $\widetilde{C}_l={\cal H}^l\cap \widetilde{C}(\Gamma)$ l'ensemble relev\'e \`a $\M$. 
Les images de ${\cal H}^l$ par $\Gamma$ sont disjointes ou confondues. 
De plus, si $l\neq l'$, les orbites $\Gamma {\cal H}^l$ et $\Gamma {\cal H}^{l'}$
sont disjointes.

En pratique, nous consid\'ererons plut\^ot la d\'ecomposition ci-dessus 
sur $\M$, not\'ee de la mani\`ere suivante:
$$
\widetilde{C}(\Gamma)=
\Gamma\,\widetilde{C}_0\sqcup \Gamma\,\widetilde{C}_1\dots\sqcup \Gamma\,\widetilde{C}_k.
$$

Notons en particulier que si $o$ est un point fix\'e de $\widetilde{C}_0$, 
son orbite $\Gamma o$ reste dans $\Gamma \widetilde{C}_0$, 
et n'intersecte donc pas les  orbites  $\Gamma{\cal H}^l$ des horoboules ${\cal H}^l$.


\subsection{Projections}


Dans toute la suite, les g\'eod\'esiques seront param\'etr\'ees \`a vitesse $1$, 
et le param\'etrage d'un rayon $[x\xi)$, avec $x\in\M$ et $\xi\in\D$ sera not\'e $(\xi_x(t))_{t\ge 0}$. 
Nous renvoyons \`a Bowditch \cite{bowditch} pour un expos\'e complet.

Rappelons que sur une vari\'et\'e d'Hadamard $\M$, i.e. une 
vari\'et\'e riemannienne simplement connexe \`a courbure n\'egative ou nulle 
(ou plus g\'en\'eralement 
sur un espace $CAT(0)$), la fonction $t\to d(x,c(t))$, 
distance d'un point $x\in\M$ \`a une g\'eod\'esique $(c(t))_{t\in\R}$ est 
une fonction propre et strictement convexe. 
Ceci permet de d\'efinir le {\em projet\'e} de $x$ sur la g\'eod\'esique 
$c$ comme l'unique point 
$c(t_0)$ qui r\'ealise le minimum de $c(t)$.

Si de plus $\M$ est une vari\'et\'e \`a courbure major\'ee par $-1$ (ou plus g\'en\'eralement un espace $CAT(-1)$),
alors ceci s'\'etend aux points $\xi$ du bord $\D$ de $\M$.
Plus pr\'ecis\'ement, si $\xi\in\D$, alors pour tout $y\in\M$ fix\'e, 
la fonction $t\to \beta_{\xi}(c(t),y)$ est strictement convexe, et l'instant $t_0$ o\`u 
elle atteint son minimum ne d\'epend pas de $y$. 
Le projet\'e $c(t_0)$ de $\xi$ sur $(c(t))_{t\in\R}$ est donc encore bien d\'efini.

Pour les m\^emes raisons de convexit\'e, on peut \'egalement d\'efinir  
le projet\'e d'un point $x\in\M\cup\D$ sur un segment g\'eod\'esique $[y z]$, 
ou un rayon g\'eod\'esique $[y \xi)$.

Remarquons enfin que si $[x \xi)$ et $[x \eta)$ sont deux rayons g\'eod\'esiques,
 alors les fonctions $t\in\R^+\to d(\xi_x(t),\eta_x(t))$ et $t\in\R^+\to d(\xi_x(t),[x\eta))$
sont \'egalement strictement convexes.

\subsection{Espaces hyperboliques au sens de Gromov}

La structure de vari\'et\'e riemannienne de $M$ n'est pas essentielle dans tout ce travail. 
En revanche, nous utiliserons de mani\`ere cruciale le fait que $\M$ est un espace hyperbolique 
au sens de Gromov, cadre d\'etaill\'e ci-dessous, dans lequel tout ce qui pr\'ec\`ede reste vrai.

Rappelons d'abord que si $(a,b,c)$ est un triangle g\'eod\'esique de $\M$, 
le {\em triangle int\'erieur} $(p,q,r)$ \`a $(a,b,c)$ est
 d\'efini comme l'unique triangle v\'erifiant $p\in[b c]$, $q\in[a c]$, $r\in[a b]$, 
et $d(a,r)=d(a,q)$, $d(b,p)=d(b,r)$ et $d(c,p)=d(c,q)$. 
Si un (ou plusieurs) des sommets, par exemple $a$ est sur le bord $\D$,
le triangle int\'erieur reste bien d\'efini, \`a condition de
remplacer la condition $d(a,q)=d(a,r)$ par $\beta_a(q,r)=0$.
Autrement dit, dans ce cas, les points $q$ et $r$ sont sur 
la m\^eme horosph\`ere centr\'ee en $a$.
Nous noterons aussi $p',q',r'$ les projet\'es respectifs 
des sommets $a,b$ et $c$ sur le c\^ot\'e oppos\'e.

\begin{figure}[!ht]
\begin{center}
\input{triangle_int.pstex_t}
\caption{Les triangles de $\M$ sont fins}
\end{center}
\end{figure}
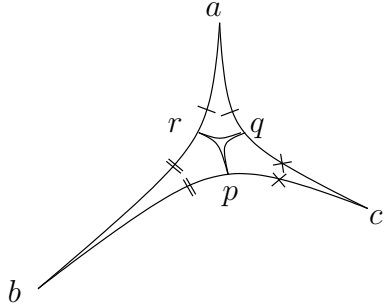

Le fait que la courbure de $\M$ soit major\'ee par $-1$ 
s'exprime dans la proposition suivante:

\begin{prop}[\cite{GH}, \cite{CDP}]\label{trianglesfins} 
Il existe une constante $\alpha\ge 0$, telle que $\M$ est
un {\em espace $\alpha$-hyperbolique au sens de Gromov}: 
tout triangle g\'eod\'esique $(a, b, c)$ de $\M\cup\D$ a un
triangle int\'erieur $(p,q,r)$ de diam\`etre inf\'erieur \`a $\alpha$.\\
De plus, on peut choisir $\alpha$ de sorte que les distances 
$d(p,p')$, $d(q,q')$ et $d(r,r')$ entre les projet\'es des sommets et les sommets du triangle int\'erieur soient toutes inf\'erieures \`a $\alpha$.
\end{prop}

L'hyperbolicit\'e permet de bien contr\^oler le d\'efaut d'\'egalit\'e 
dans l'in\'egalit\'e triangulaire $d(b,c)\le d(a,b)+d(a,c)$.
D'un point de vue riemannien, si l'angle au sommet $a$ est minor\'e par $\theta>0$, 
il existe une constante $C(\theta)$, telle que 
$d(a,b)+d(a,c)-d(b,c)\le C(\theta)$. 
R\'eciproquement, si $d(a,b)+d(a,c)-d(b,c)\le c$, 
et si les c\^ot\'es du triangle ne sont pas trop petits, 
alors l'angle au sommet $a$ est minor\'e par une constante $\theta(c)$. 
Bien que plus visuels, les angles riemanniens sont en fait plus difficiles 
\`a manipuler que des distances, et nous pr\'ef\'ererons donc l'\'enonc\'e ci-dessous: 

\begin{lem}\label{angleausommetbis}
Soit $(a,b,c)$ un triangle g\'eod\'esique de $\M$, et $C>0$. 
Si $d(a,[bc])\le C$, alors 
$\displaystyle d(a,b)+d(a,c)-d(b,c)\le 2C$. 
R\'eciproquement, 
si $\displaystyle d(a,b)+d(a,c)-d(b,c)\le C$, 
alors $d(a,[bc])\le \frac{C}{2}+\alpha$.\\
Ceci reste vrai si $b$ ou $c$ est un point du bord, en rempla\c cant $d(a,b)+d(a,c)-d(b,c)$ par 
$\beta_b(a,x)+\beta_c(a,x)$, pour tout $x\in (b,c)$.
\end{lem}

\begin{demo}
Par d\'efinition du triangle int\'erieur, on a 
$d(a,b)+d(a,c)-d(b,c)=2d(a,r)$. Puisque
$d(a,[bc])\le d(a,r)+\alpha$, on en d\'eduit la deuxi\`eme partie du lemme.\\
Notons maintenant $p'$ le projet\'e de $a$ sur $[bc]$, et supposons que $d(a,p')=d(a,[bc])\le C$.
Alors on  a:
$$
0\le d(a,b)+d(a,c)-d(b,c)=d(a,b)-d(p',b)+d(a,c)-d(p',c)\le 2d(a,p')\le 2C.
$$
\end{demo}

Si le projet\'e de $a$ sur $[bc]$ est $b$, 
alors il est clair que l'angle au sommet $b$ est sup\'erieur ou \'egal \`a $\pi/2$. 
En termes de distances, ceci se r\'e\'ecrit:

\begin{lem}\label{projetebis}
Soit $(a,b,c)$ un triangle g\'eod\'esique de $\M$, et $p'$ le projet\'e de $a$ sur le c\^ot\'e  $[bc]$. 
Si $p'=b$, alors la distance de $b$ \`a $[ac]$ est inf\'erieure \`a $2\alpha$.
R\'eciproquement, si cette distance est inf\'erieure \`a $2\alpha$, 
alors la distance de $p'$ \`a $b$ est au plus $3\alpha$.
\end{lem}

Ce lemme se g\'en\'eralise imm\'ediatement au cas o\`u $a$ est un point du bord, 
\`a condition d'utiliser les fonctions de Busemann au lieu de la distance.

\begin{demo}
D'apr\`es la proposition \ref{trianglesfins}, 
on a d'une part $d(b,[ac])\le d(b,p)+\alpha$, et d'autre part 
$d(p,p')\le \alpha$. Si $b=p'$, on trouve imm\'ediatement $d(b,[ac])\le 2\alpha$.

R\'eciproquement, si $d(b,[ac])\le 2\alpha$, le lemme \ref{angleausommetbis}
donne $2d(b,p)=d(b,a)+d(b,c)-d(a,c)\le 4\alpha$, et 
d'apr\`es la proposition \ref{trianglesfins}, $d(p,p')\le \alpha$, ce qui
montre bien que la distance de $b$ \`a $p'$ est au plus $3\alpha$.
\end{demo}

\subsection{Voisinages d'un point du bord}\label{Voisinagesdunpointdubord}

Il existe plusieurs familles \'equivalentes de voisinages d'un point $\xi$ du bord. 
Nous n'en utiliserons qu'une, mais nous les d\'efinissons toutes, de mani\`ere \`a 
donner une bonne image de ce que signifie << \^etre proche >> pour deux points du bord.
Dans ce qui suit, on consid\`ere un point $x$ de $\M$, 
un point $\xi$ du bord $\D$, et on note toujours 
$(\xi_x(t))_{t\ge 0}$ la param\'etrisation \`a vitesse $1$ du rayon $[x \xi)$.\\

Les voisinages les plus c\'el\`ebres  de $\xi$ sont sans doute les {\em ombres} de Sullivan. On fixe un 
r\'eel $r>0$, et pour $x\in\M$, $\xi\in\D$ et $t>0$, on d\'efinit 
${\cal O}(x,\xi,t)$ comme l'ombre faite sur le bord 
par la boule de centre $\xi_x(t)$ et de rayon $r$ vue du point $x$, 
soit encore l'ensemble des $\eta\in\D$ tels que le rayon $[x \eta)$ intersecte la boule $B(\xi_x(t),r)$.\\

On peut \'egalement consid\'erer les boules de la famille $(d_x)_{x\in\M}$ des {\em distances de Gromov}, 
ou distances visuelles sur le bord, 
d\'efinies pour tout $x\in\M$ et tous $(\xi,\eta)\in\DD$ par:
$$
d_x(\xi,\eta)=
\exp\left(-\frac{1}{2}\beta_{\xi}(x,y)-\frac{1}{2}\beta_{\eta}(x,y)\right),\quad\mbox{avec}\; y\in(\xi\eta)\,.
$$
C'est ici que nous nous servons du fait que la courbure de $\M$ est major\'ee par $-1$; 
en effet, si elle est inf\'erieure \`a $-b^2$ avec $b< 1$, 
les quantit\'es ci-dessus sont toujours d\'efinies 
mais ne satisfont plus l'in\'egalit\'e triangulaire, voir \cite{bourdon}.
Nous noterons $B_x(\xi,r)$ la boule de centre $\xi$ et de rayon $r$ pour la distance $d_x$.

Les voisinages de Hamenst\"adt 
sont d\'efinis comme suit: 
${\cal D}(x,\xi,t)$ est l'ensemble des points $\eta\in\D$ tels que la 
distance de $\xi_x(t)$ \`a $\eta_x(t)$ est inf\'erieure \`a $\alpha$.
La proposition suivante montre que ces trois familles de voisinages sont 
pratiquement les m\^emes:

\begin{prop}[Kaimanovich, \,\cite{kaim2}]\label{comp}Soit $\M$ une vari\'et\'e de Hadamard \`a courbures sectionnelles major\'ees par $-1$.
Il existe des constantes 
$c_1\ge 1$ et $c_2\ge 1$ telles que pour tous $x\in\M$, $\xi\in\D$ et $t\ge 0$, on ait: 
$$
B_x(\xi,\frac{1}{c_1}e^{-t})\subset {\cal O}(x,\xi,t)\subset B_x(\xi,c_1\,e^{-t}),\quad \mbox{et}
$$
$$
B_x(\xi,\frac{1}{c_2}e^{-t})\subset {\cal D}(x,\xi,t)\subset B_x(\xi,c_2\,e^{-t})\,.
$$
\end{prop}

Pour les besoins de la preuve du th\'eor\`eme \ref{courbvariable}, il sera utile de travailler avec un dernier type de voisinages.
D\'efinissons $V(x,\xi,0)$ comme l'ensemble des points $\eta\in\D$ dont le 
projet\'e sur la g\'eod\'esique $(x \xi)$ appartient en fait au rayon $[x \xi)$, 
et si $t\ge 0$, $V(x,\xi,t)\subset V(x,\xi,0)$ comme l'ensemble de ceux qui se projettent 
sur $[\xi_x(t) \,\xi)$, i.e. \`a distance strictement sup\'erieure \`a
 $t$ de $x$. Le lemme suivant montre qu'ils sont encore comparables aux voisinages d\'efinis ci-dessus:

\begin{lem}\label{comparables}
Pour tout $x\in\M$, $\xi\in\D$, et $t\ge 2\alpha$, on a:
$$
V(x,\xi,t+\alpha)\subset {\cal D}(x,\xi,t)\subset V(x,\xi,t-2\alpha).
$$
\end{lem}

Notons que m\^eme en courbure constante \'egale \`a $-1$, 
ces voisinages ne co\"incident pas exactement les uns avec les autres. 
Par exemple, on peut montrer que si $t\ge 0$, 
$\displaystyle V(x,\xi,t)=B_x(\xi,\frac{e^{-t}}{\sqrt{1+e^{-2t}}})$.

\begin{demon}{du lemme \ref{comparables}}
Soit tout d'abord $\eta$ un point de ${\cal D}(x,\xi,t)$, qui v\'erifie donc $d(\xi_x(t),\eta_x(t))\le~\alpha$.
Notons $u_0$ le minimum de la fonction  $\varphi(u)=\beta_{\eta}(\xi_x(u),\eta_x(t))$. 
Le point  $\xi_x(u_0)$ est donc le projet\'e de $\eta$ sur le rayon $[x\xi)$.
L'in\'egalit\'e trianglaire donne $\varphi(t)=\beta_{\eta}(\xi_x(t),\eta_x(t))\le d(\xi_x(t),\eta_x(t))\le \alpha$.
D'autre part, si $u\le t$, on a 
$$
\varphi(u)=\beta_{\eta}(\xi_x(u),\eta_x(t))=\beta_{\eta}(\xi_x(u),\eta_x(u))+\beta_{\eta}(\eta_x(u),\eta_x(t)) \ge -\alpha+t-u.
$$  
Si $u<t-2\alpha$, on en d\'eduit que $\varphi(u)>\varphi(t)$. 
Par stricte convexit\'e de $\varphi$, le minimum de $\varphi$ est n\'ecessairement atteint en $u_0\ge t-2\alpha$, 
ce qui signifie exactement que $\eta$ appartient \`a $V(x,\xi,t-2\alpha)$.

R\'eciproquement, si $\eta\in V(x,\xi,t+\alpha)$, 
d'apr\`es la proposition \ref{trianglesfins}, 
le sommet $\xi_x(s_0)$ du triangle int\'erieur \`a $(x,\xi,\eta)$ v\'erifie $s_0\ge t$, 
d'o\`u on d\'eduit\\ $d(\xi_x(t),\eta_x(t))\le d(\xi_x(s_0),\eta_x(s_0))\le \alpha$.
\end{demon}

Nous d\'etaillons maintenant certaines propri\'et\'es de ces ensembles 
$V(x,\xi,t)$ qui nous serviront au cours de la preuve du
th\'eor\`eme \ref{courbvariable}.
La premi\`ere est \'el\'ementaire:

\begin{lem}\label{anglebis}
Pour tout $x\in\M$, $\xi\in\D$ et $t\ge 0$, si $\eta\in V(x,\xi,t)$, 
alors la distance de $\xi_x(t)$ \`a $(x\eta)$ est inf\'erieure \`a $2\alpha$. 
En particulier, $\displaystyle  t-4\alpha\le\beta_{\eta}(x,\xi_x(t))\le t$.
\end{lem}

\begin{demo}L'in\'egalit\'e triangulaire donne $\beta_{\eta}(x,\xi_x(t))\le t$. 
Le lemme \ref{projetebis} implique $d(\xi_x(t),[x\eta))\le 2\alpha$.
Par le lemme \ref{angleausommetbis}, nous en d\'eduisons $\beta_{\eta}(x,\xi_x(t))\ge t-4\alpha$.
\end{demo}

\begin{figure}[ht!]
\begin{center}
\input{voisinage.pstex_t}
\caption{$V(x,\xi,t)$}
\end{center}
\end{figure}
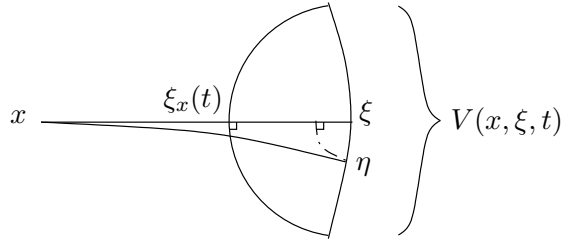

Dans le lemme suivant, nous montrons comment varient ces ensembles $V(x,\xi,t)$ quand on 
fait varier l\'eg\`erement $\xi$ ou $x$.  

Notons le fait \'evident suivant: si $\eta\in B_x(\xi,r/2)$, 
alors $\displaystyle B_x(\eta,r/2)\subset B_x(\xi,r)\subset B_x(\eta,3r/2)$.
La partie b- du lemme ci-dessous est une simple reformulation de cette inclusion, 
et en d\'ecoule directement si on utilise la proposition \ref{comp} et le lemme \ref{comparables}. 

\begin{lem}\label{encadrement}  
{\bf a-} Soient $x\in\M$, $\xi\in\D$, et $t\ge 6\alpha$. Alors pour tout $\eta\in {\cal D}(x,\xi,t)$ 
(i.e. tel que $d(\xi_x(t),\eta_x(t))\le\alpha$), on a 
$$
V(x,\eta,t)\subset V(x,\xi,t-6\alpha).
$$
{\bf b-}
Posons  $K_1=6\alpha>0$. Soient $x\in\M$,  $\xi\in\D$ et $t\ge K_1$. Pour tout $\eta\in V(x,\xi, t+K_1+\alpha)$, on a:
$$
V(x,\eta,t+K_1)\subset V(x,\xi,t)\subset V(x,\eta, t-K_1)
$$
{\bf c- }Pour tout $D>0$, notons $K_2=K_2(D)=2D+4\alpha>0$. Soient 
$x$ et $y$ des points de $\M$ tels que $d(x,y)\le D$,  $\xi\in\D$, et $t\ge K_2$. Alors:
$$
V(x,\xi,t+K_2)\subset V(y,\xi,t)\subset V(x,\xi,t-K_2)\,.
$$
\end{lem}

\begin{demo}{\bf a- }Soient $\eta\in {\cal D}(x,\xi,t)$ et $\zeta\in V(x,\eta,t)$.
Notons $\xi_x(u_0)$ le projet\'e de $\zeta$ sur le rayon $[x\xi)$.
Par d\'efinition, $u_0$ minimise la fonction 
$u\to \varphi(u)=\beta_{\zeta}(\xi_x(u),\eta_x(t))$. 
Or 
$$
\varphi(t)=\beta_{\zeta}(\xi_x(t),\eta_x(t))\le d(\xi_x(t),\eta_x(t))\le \alpha 
$$ 
par hypoth\`ese. 
Par ailleurs, si $u\le t$, l'in\'egalit\'e triangulaire et le lemme \ref{anglebis} donnent:
$$
\varphi(u)=\beta_{\zeta}(\xi_x(u),\eta_x(u))+\beta_{\zeta}(\eta_x(u),\eta_x(t))\ge -\alpha+t-u-4\alpha.
$$
Autrement dit, si $0\le u<t-6\alpha$, $\varphi(u)>\varphi(t)$. 
La stricte convexit\'e de $\varphi$ implique alors que son minimum 
est atteint en $u_0\ge t-6\alpha$. Ceci signifie exactement que $\zeta\in V(x,\xi,t-6\alpha)$.

\noindent
{\bf b- }Si $\eta\in V(x,\xi, t+K_1+\alpha)$, le lemme \ref{comparables} 
donne 
$$d(\xi_x(t),\eta_x(t))\le d(\xi_x(t+K_1),\eta_x(t+K_1))\le \alpha.$$ 
On en d\'eduit, en appliquant deux fois le {\bf a}, que 
$V(x,\xi,t)\subset~V(x,\eta, t-K_1)$ et $V(x,\eta,t+K_1)\subset~V(x,\xi,t)$, 
ce qui est le r\'esultat souhait\'e.

\noindent
{\bf c- }Soit $\zeta\in V(y,\xi,t)$, avec $t\ge K_2$.
 Posons $\varphi(u)=\beta_{\zeta}(\xi_x(u),\xi_y(t))$.
Si $d(x,y)\le D$, comme les rayons $[x\xi)$ et $[y\xi)$ sont asymptotes, on a 
$d(\xi_x(u),\xi_y(u))\le d(x,y)\le D$ pour tout $u\ge 0$. 
On en d\'eduit $|\beta_{\zeta}(\xi_x(t),\xi_y(t))|\le D$. 
Et pour tout $0\le u\le t$, d'apr\`es le lemme \ref{anglebis}, 
$$
\varphi(u)=\beta_{\zeta}(\xi_x(u),\xi_y(u))+\beta_{\zeta}(\xi_y(u),\xi_y(t))\ge -D+t-u-4\alpha.
$$ 
Comme $\varphi$ est strictement convexe, le projet\'e 
$\xi_x(u_0)$ de $\zeta$ sur $(x\xi)$ v\'erifie $u_0\ge t-2D-4\alpha=t-K_2$, 
ce qui montre que $V(y,\xi,t)\subset V(x,\xi,t-K_2)$.
L'autre inclusion se prouve de la m\^eme mani\`ere.
\end{demo}

Pour finir ce paragraphe, \'enon\c cons un dernier lemme qui 
nous servira par la suite, et qui d\'ecoule imm\'ediatement des 
lemmes \ref{encadrement} c, \ref{comparables} et de la proposition \ref{comp}:

\begin{lem}\label{epsilon}
Soit $D>0$ fix\'e. 
Il existe $\varepsilon>0$, tel que pour tout $o\in\widetilde{M}$, et 
tout $x\in B(o,D)$ et tout $\xi\in\D$:
$$V(x, \xi, 0)\supset B_{o}(\xi,\varepsilon)$$
\end{lem}



\subsection{Action d'une isom\'etrie parabolique sur le bord}\label{s25}


Nous allons maintenant nous servir de ces ensembles $V(x,\xi,t)$ pour comprendre 
l'action d'une isom\'etrie parabolique de $\M$ sur le bord priv\'e de son point fixe,
ainsi que sur les horosph\`eres qu'elle stabilise (i.e. celles qui sont centr\'ees en son point fixe).
Une isom\'etrie parabolique, quand elle est it\'er\'ee,
 attire tous les points du bord d'une part, et de $\M$ d'autre part vers son point fixe. 
Ce que dit le lemme ci-dessous, c'est qu'on peut quantifier et relier entre eux ces d\'eplacements.

\begin{lem}\label{lemme29}
Soient $o\in\M$, $\xi\in\D$, et $K$ un compact de $\D\setminus\{\xi\}$. 
Il existe une constante $K_3>0$, telle que pour toute isom\'etrie parabolique $p$ fixant $\xi$, et tout $t\ge K_3$, 
on a:\\
{\bf a- }Si $d(o,po)\ge 2t$, alors $pK\subset V(o,\xi,t-K_3)$ et pour tout $\eta\in K$, 
$$|\beta_{p\eta}(\xi(t),p\xi(t))-d(o,po)+2t|\le 2K_3.$$
{\bf b- } Si $d(o,po)\le 2t$, alors $p K\cap V(o,\xi,t+K_3)=\emptyset$ et 
$|\beta_{p\eta}(\xi(t),p\xi(t))|\le 2K_3$ pour tout $\eta\in K$.
\end{lem}

\begin{figure}[ht!]
\begin{center}
\input{lemme29.pstex_t}
\caption{Lemme \ref{lemme29}}
\end{center}
\end{figure}

\begin{demo}
Si $\eta\in K$, notons $y_{\eta}$ l'intersection de la g\'eod\'esique $(\xi\eta)$ 
et de l'horosph\`ere centr\'ee en $\xi$ et passant par $o$. 
Par compacit\'e de $K$, on a:
$$
C(K,o,\xi):=\sup_{\eta\in K}\,d(o,y_{\eta})<+\infty.
$$
Notons $(\xi(s))_{s\ge 0}$ le rayon $[o\xi)$
et $(\xi_{\eta}(s))_{s\ge 0\in\R}$ la param\'etrisation de la 
g\'eod\'esique $(\eta\xi)$ positive sur le rayon $[y_{\eta}\,\xi)$.
On a donc $y_{\eta}=\xi_{\eta}(0)$ et $o=\xi(0)$.

Consid\'erons le triangle $(o,\xi,p\eta)$. 
Comme $\beta_{\xi}(o,py_{\eta})=0$, les sommets des c\^ot\'es 
$[o\xi)$ et $[p\eta\,\xi)$ de son triangle int\'erieur s'\'ecrivent 
$\xi(s_0)$ et $\xi_{p\eta}(s_0)=p\xi_{\eta}(s_0)$, avec $s_0\ge 0$. 
L'in\'egalit\'e triangulaire et la d\'efinition du triangle int\'erieur 
donnent facilement $\displaystyle 2s_0-\alpha\le d(o,py_{\eta})\le 2s_0+\alpha$. 
D'o\`u on d\'eduit, puisque $d(po,py_{\eta})=d(o,y_{\eta})\le C(K,o,\xi)$,
\begin{eqnarray}\label{encadrementdistance}
2s_0-C(K,o,\xi)-\alpha\le d(o,po)\le 2s_0+C(K,o,\xi)+\alpha
\end{eqnarray}
D'apr\`es la proposition \ref{trianglesfins}, 
le projet\'e de $p\eta$ sur $[o\xi)$ est \`a distance au plus $\alpha$ de $\xi(s_0)$. 
Posons $K_3=C(K,o,\xi)/2+3\alpha/2$.
Avec l'encadrement (\ref{encadrementdistance}), l'in\'egalit\'e
$d(o,po)\ge 2t$ implique $p\eta\in V(o,\xi,t-K_3)$. 
Ceci \'etant vrai pour tout $\eta\in K$, on a bien $pK\subset  V(o,\xi,t-K_3)$.\\
Si $d(o,po)\le 2t$, on montre de m\^eme que $pK\cap V(o,\xi,t+K_3)=\emptyset$.\\

Il nous reste maintenant \`a estimer la quantit\'e $\beta_{p\eta}(\xi(t),p\xi(t))$ pour $\eta\in\D$.
Remarquons d'abord que 
$$
|\beta_{p\eta}(p\xi_{\eta}(t),p\xi(t))|\le d(\xi_{\eta}(t),\xi(t))\le d(y_{\eta},o)\le C(K,o,\xi).
$$
Il nous suffit donc d'estimer $\beta_{p\eta}(\xi(t),p\xi_{\eta}(t))$.
Si $s_0\le t$, alors 
$$
|\beta_{p\eta}(\xi(t),p\xi_{\eta}(t))|\le d(\xi(t),p\xi_{\eta}(t))\le d(\xi(s_0),p\xi_{\eta}(s_0))\le\alpha.
$$

Supposons maintenant  $s_0\ge t$, et notons $(c(t))_{t\ge 0}$ la param\'etrisation du rayon g\'eod\'esique 
$[o\,p\eta)$. Le troisi\`eme sommet du triangle int\'erieur \`a  $(o,\xi, p\eta)$ est alors $c(s_0)$.
Par d\'efinition de ce triangle int\'erieur, on a alors $\beta_{p\eta}(c(s_0),p\xi_{\eta}(s_0))=0$, 
d'o\`u on d\'eduit que $\beta_{p\eta}(c(2s_0-t),p\xi_{\eta}(t))=0$.
D'autre part, on a aussi $|\beta_{p\eta}(\xi(t),c(t))|\le d(\xi(t),c(t))\le\alpha$.
A l'aide de toutes ces in\'egalit\'es, on obtient:
$$
\beta_{p\eta}(\xi(t),p\xi_{\eta}(t))=
\beta_{p\eta}(\xi(t),c(t))+\beta_{p\eta}(c(t),c(2s_0-t))=\beta_{p\eta}(\xi(t),c(t))+2(s_0-t), 
$$
puis $|\beta_{p\eta}(\xi(t),p\xi(t))-(2s_0-2t)|\le \alpha+C(K,o,\xi)$.
Finalement, en utilisant la notation $u^+=\max(u,0)$ pour $u\in\R$, 
on peut rassembler les deux cas $s_0\le t$ et $s_0\ge t$ ci-dessus pour obtenir pour tout $t\ge 0$:
\begin{eqnarray}\label{encadrementBusemann}
(2s_0-2t)^+-\alpha-C(K,o,\xi)\le \beta_{p\eta}(\xi(t),p\xi(t))\le (2s_0-2t)^+ +\alpha+C(K,o,\xi).
\end{eqnarray}

\noindent
Une manipulation \'el\'ementaire des encadrements (\ref{encadrementdistance}) et (\ref{encadrementBusemann})
permet de conclure la preuve du lemme.
\end{demo}


\section{Mesures conformes et lemme de l'ombre}\label{s3}

\subsection{Densit\'es conformes}

Soit \`a pr\'esent $\Gamma$ un groupe discret d'isom\'etries de $\M$.
Une {\em densit\'e conforme} $\Gamma$-invariante de dimension $\delta>0$ sur $\D$ est 
une famille $\nu=(\nu_x)_{x\in\M}$ de mesures finies sur $\D$ qui 
v\'erifient la condition d'invariance $\nu_{\gamma x}=\gamma_{*}\nu_x$ 
pour tout $\gamma\in\Gamma$, et pour tous $(x,y)\in\M^2$ et $\xi\in\D$:
$$
\frac{d\nu_x}{d\nu_y}(\xi)=\exp(-\delta\beta_{\xi}(x,y))\,,\quad \nu_y - \mbox{p.s.} 
$$
Soit $o\in\M$ un point fix\'e dans toute la suite de ce paragraphe.
La densit\'e $\nu$ sera dite {\em normalis\'ee} si $\nu_o$ est une probabilit\'e.

L'{\em exposant critique} d'un groupe discret d'isom\'etries $\Gamma$ de $\M$ est d\'efini par:
$$
\delta_{\Gamma}=\limsup_{T\to\infty}\,\frac{1}{T}\,\log\,\sharp\{\gamma\in\Gamma, d(o,\gamma o)\le T\}\,.
$$
C'est encore l'exposant critique de la s\'erie de Poincar\'e de $\Gamma$:
$$
P(\Gamma,o,s)=\sum_{\gamma \in \Gamma}\,e^{-sd(o,\gamma o)}
$$
Le groupe $\Gamma$ est dit {\em de type divergent} si la s\'erie ci-dessus diverge en $s=\delta_{\Gamma}$. 
Si $\Gamma$ est {\it non \'el\'ementaire},  c'est-\`a-dire si 
$\sharp\Lambda_{\Gamma}=+\infty$, alors $\delta_{\Gamma}$ est strictement positif, et 
Roblin \cite{Roblin3} a d\'emontr\'e que la limite sup\'erieure ci-dessus est en fait une limite.
Mais un tel r\'esultat est inconnu pour les groupes \'el\'ementaires, 
par exemple les sous-groupes paraboliques de $\Gamma$. 
Pour d\'emontrer le th\'eor\`eme \ref{courbvariable}, 
nous aurons besoin justement de 
supposer que c'est bien une limite, et plus encore: 
nous supposerons que pour tout sous-groupe parabolique $\Pi$ de $\Gamma$, 
il existe une constante $D\ge 1$, telle que 
\begin{eqnarray}\label{hypothesecroissance}
\frac{1}{D}\,\exp(\delta_{\Pi}T)\,\le\,\sharp\{p\in\Pi, d(o,p o)\in [ T,T+1[\}\,\le D\,\exp(\delta_{\Pi}T).
\end{eqnarray}

Cette hypoth\`ese est apparue \`a plusieurs reprises dans la litt\'erature, d'abord dans \cite{Paulin},
 mais aussi \cite{Franchi} par exemple, et semble n\'ecessaire 
\`a chaque fois qu'on souhaite obtenir des estim\'ees pr\'ecises \`a l'int\'erieur des cusps. 
La proposition ci-dessous montre qu'elle est v\'erifi\'ee d\`es que les cusps sont << raisonnables >>.

\begin{prop}\label{satisfaite}
L'hypoth\`ese (\ref{hypothesecroissance}) est v\'erifi\'ee d\`es que 
les cusps de $M$ sont isom\'etriques aux cusps d'une vari\'et\'e localement sym\'etrique de rang un.
\end{prop}

\begin{demo}
Pour chaque sous-groupe parabolique $\Pi$, 
on peut v\'erifier cette hypoth\`ese en supposant que $o$ 
est dans le cusp consid\'er\'e. 
Or le sous-groupe parabolique $\Pi$ stabilise le cusp, 
donc $\Pi o$ reste dans les relev\'es du cusp. 
La preuve consiste alors simplement \`a v\'erifier que 
dans les espaces localement sym\'etriques, 
cette condition de croissance est v\'erifi\'ee. 
Dans l'espace hyperbolique r\'eel, 
c'est un simple calcul utilisant d'une part le fait qu'une isom\'etrie parabolique $p$ 
fixant le point $\infty$ (dans le mod\`ele du demi-espace sup\'erieur) agit 
par translation euclidienne sur les horosph\`eres horizontales, et d'autre part l'\'equivalent 
$d(o,po)\simeq 2\log d_{\rm eucl}(o,po)$, o\`u $o$ est le point \`a hauteur (euclidienne) $1$ sur l'axe vertical.
Dans les espaces hyperboliques exotiques, il est encore possible d'obtenir un \'equivalent exact de cette distance $d(o,po)$ 
(voir Corlette-Iozzi \cite[Formule 3.5]{CI}), ce qui permet de conclure (voir aussi une preuve compl\`ete dans  \cite[lemma 3.5]{newberger}).  
\end{demo}

Cette hypoth\`ese implique de mani\`ere imm\'ediate 
que tout sous-groupe parabolique $\Pi$ de $\Gamma$ est divergent. 
Des travaux de Dal'bo, Otal et Peign\'e \cite{DOP}, on d\'eduit le r\'esultat suivant:

\begin{theo}[Dal'bo-Otal-Peign\'e \cite{DOP}]Si $\Gamma$ est un groupe g\'eom\'etriquement fini 
dont tout sous-groupe parabolique $\Pi$ est divergent, 
alors $\delta_{\Pi}<\delta_{\Gamma}$ et le groupe $\Gamma$ est lui-m\^eme divergent.
\end{theo}

Si $\Gamma$ est un groupe discret non \'el\'ementaire, 
un r\'esultat originellement d\^u \`a Patterson \cite{Patterson}
assure l'existence d'une densit\'e conforme invariante de dimension 
$\delta_{\Gamma}$ \`a support $\Lambda_{\Gamma}$. 
L'\'etude d'abord faite par Sullivan \cite{Sull} \cite{Sullivan} en courbure constante, 
puis g\'en\'eralis\'ee en courbure variable\cite{yue} montre aussi que
toute densit\'e $\delta$-conforme invariante v\'erifie $\delta\ge\delta_{\Gamma}$. 
Lorsque $\Gamma$ est divergent, ce qui est v\'erifi\'e d'apr\`es le 
th\'eor\`eme ci-dessus sous notre hypoth\`ese (\ref{hypothesecroissance}), on a de plus:

\begin{prop}\label{uniquedensite} 
Si $\Gamma$ est non \'el\'ementaire divergent, il existe une unique densit\'e 
conforme invariante normalis\'ee $\nu=(\nu_x)_{x\in\M}$
 de dimension $\delta_{\Gamma}$ \`a support $\Lambda_{\Gamma}$, et elle 
 est sans atomes. 
La mesure $\nu_o$ est alors appel\'ee {\em mesure de Patterson}.
\end{prop}

Dans le paragraphe ci-dessous, nous travaillerons avec une densit\'e conforme 
invariante de dimension quelconque, mais la densit\'e donn\'ee par la proposition 
ci-dessus sera privil\'egi\'ee dans tout le paragraphe \ref{s4}.


\subsection{Lemme de l'ombre}\label{Lemmedelombre}


Dans le cas d'une vari\'et\'e compacte \`a courbure n\'egative, le lemme de l'ombre, 
originellement d\^u \`a Sullivan \cite{Sull} en courbure constante, 
permet de montrer que pour tout $x\in\M$,
la mesure de Patterson $\nu_x$ de $\Gamma$ est la mesure de Hausdorff 
(de dimension $\delta_{\Gamma}$) pour la distance de Gromov $d_x$ sur $\D$.
En revanche, dans le cas d'une vari\'et\'e non compacte, 
la description de $\nu$ est beaucoup plus d\'elicate.
Le cas des vari\'et\'es g\'eom\'etriquement finies de courbure constante 
fut d\'ecrit par Stratmann et Velani \cite{SV}, et r\'ecemment
g\'en\'eralis\'e aux vari\'et\'es localement sym\'etriques de rang $1$ par F. Newberger dans \cite{newberger}.
Nous montrerons que cette description est encore valide pour 
les vari\'et\'es g\'eom\'etriquement finies de courbure variable 
sous la condition (\ref{hypothesecroissance}).
Pour cela, nous adapterons la preuve de M. Peign\'e \cite{peigne} 
du r\'esultat de Stratmann et Velani.

Les estim\'ees du th\'eor\`eme ci-dessous  se d\'emontrent cusp par cusp, 
nous supposerons donc dans ce paragraphe que $M$ n'a qu'un seul cusp not\'e $C_1$. \\

Rappelons les notations du paragraphe \ref{s21}.
Le coeur de Nielsen $N_{\Gamma}$ est un convexe qui 
se d\'ecompose en $N_{\Gamma}=C_0\sqcup C_1$, 
o\`u $C_0$ est la partie compacte de la vari\'et\'e, 
de diam\`etre inf\'erieur \`a $\Delta$, 
et le cusp $C_1$ est le quotient de $\tilde{C}_1:={\cal H}_1\cap \widetilde{C}(\Gamma)$ 
par $\Pi$, o\`u ${\cal H}_1$ est une horoboule centr\'ee au 
point parabolique $\xi_1\in\Lambda_{\Gamma}$, 
et $\Pi$ est le stabilisateur de $\xi_1$ dans $\Gamma$.
On peut alors choisir un relev\'e  $\widetilde{C}_0$ de $C_0$ \`a $\M$, 
qui soit connexe et de diam\`etre inf\'erieur \`a $\Delta$, 
 et tel que les adh\'erences de $\widetilde{C}_0$ et $\widetilde{C}_1$ s'intersectent.
On choisira une famille $\{\gamma_i\}_{i\ge 0}$ de repr\'esentants de $\Gamma/\Pi$, 
et on notera $\xi_i=\gamma_i\xi_1$ et ${\cal H}_i$ l'horoboule $\gamma_i{\cal H}_1$. On a alors 
On a alors:
$$
\Gamma \widetilde{C}_1\,=\,\sqcup_{i=1}^{+\infty}\gamma_i\widetilde{C}_1\,
=\,\sqcup_{i=1}^{+\infty}{\cal H}_i\cap \widetilde{C}(\Gamma)\,.
$$
D'autre part,  
quitte \`a faire agir le stabilisateur 
$\gamma_i\Pi\gamma_i^{-1}$ de ${\cal H}_i$, 
nous pouvons supposer que pour tout $i\in\N$, 
$\gamma_i\widetilde{C}_0$ intersecte le rayon $[o\xi_i)$ 
(et le bord $\partial {\cal H}_i$ \'evidemment).

Rappelons pour finir que si $\xi\in\D$, 
$\xi(t)$ d\'esigne toujours le point \`a distance $t$ 
de $o$ sur le rayon g\'eod\'esique $[o\xi)$.

Le th\'eor\`eme ci-dessous \'etend le Lemme de l'Ombre en courbure variable, 
permettant ainsi d'obtenir des estim\'ees pr\'ecises de la mesure des ombres sur le bord.

\begin{theo}\label{courbvariable}
Soit $M=\M/\Gamma$ une vari\'et\'e g\'eom\'etriquement finie \`a un cusp $C_1$, 
et $\Pi$ un sous-groupe parabolique de $\Gamma$, tel qu'il
existe une constante $D\ge 1$, telle que 
\begin{eqnarray*}
\frac{1}{D}\exp(\delta_{\Pi}T)\,\le\,\sharp \{p\in\Pi,d(o, po)\in[T,T+1[\}\,\le\,D\exp(\delta_{\Pi}T).
\end{eqnarray*}
Alors pour toute  densit\'e conforme $\nu=(\nu_x)_{x\in\M}$ 
$\Gamma$-invariante sans atomes normalis\'ee 
de dimension $\delta\ge \delta_{\Gamma}$ et de support 
$\Lambda_{\Gamma}$, il existe des constantes $A_0>0$ et $A_1>0$ telles que pour tout $\xi\in\Lambda_{\Gamma}$ et $t\ge 0$:\\
{\bf a- } si $\xi(t)\in\Gamma\widetilde{C}_0$, alors
$\displaystyle\quad\frac{1}{A_0}\exp(-\delta t)\le\nu_{o}(V(o,\xi,t))\,\le A_0 \exp(-\delta t)$,\\
{\bf b- } si $\xi(t)\in\Gamma \widetilde{C}_1$, alors:
$$
\frac{1}{A_1}\,e^{-\delta t+(2\delta_{\Pi}-\delta)\,d(\xi(t),\Gamma o)}\,\le\,
\nu_{o}(V(o,\xi,t))\,\le \,A_1 \,e^{-\delta t+(2\delta_{\Pi}-\delta)\,d(\xi(t),\Gamma o)} .
$$
\end{theo}

En particulier, le th\'eor\`eme s'applique \`a l'unique 
densit\'e $\delta_{\Gamma}$-conforme invariante donn\'ee 
par la proposition \ref{uniquedensite}, qui est sans atomes.


\subsection{Preuve du th\'eor\`eme \ref{courbvariable}}

Le plan de la preuve est celui de Peign\'e \cite{peigne}; 
la diff\'erence essentielle r\'eside dans la  proposition \ref{cle}, 
o\`u se manifeste la courbure variable, 
et o\`u l'hypoth\`ese (\ref{hypothesecroissance}) est utilis\'ee.
Nous rappelons toutefois toute la d\'emonstration pour la commodit\'e du lecteur.
Elle se fait en plusieurs \'etapes: 
le lemme \ref{ombreclassique} est le lemme de l'ombre classique, 
il traite le cas o\`u 
$\xi(t)$ appartient au relev\'e $\Gamma\widetilde{C}_0$ de la partie compacte.
Le r\'esultat cl\'e est la proposition \ref{cle}. 
De cette proposition d\'ecoulent le corollaire \ref{parabolique}, 
qui traite le cas o\`u $\xi=\xi_i$ est un point parabolique et $\xi(t)$ appartient \`a 
l'horoboule ${\cal H}_i$ centr\'ee en $\xi_i$, et le lemme \ref{derniercas}
qui permet de conclure la preuve du th\'eor\`eme \ref{courbvariable} dans tous les autres cas.\\

Commen\c cons par une remarque importante: 
si $\eta\in V(o,\xi,t)$, on a 
$|\beta_{\eta}(o,\xi(t))-t|\le 4\alpha$ par le lemme \ref{anglebis}. 
En utilisant la relation de conformit\'e de $\nu$, on obtient:
\begin{eqnarray}\label{simplification1}
\exp(-4\delta \alpha+\delta t)\,\le \,\frac{\nu_{\xi(t)}(V(o,\xi,t))}{\nu_{o}(V(o,\xi,t))}\,\le\, \exp(4\delta \alpha+\delta t).
\end{eqnarray}
Autrement dit, pour prouver le th\'eor\`eme \ref{courbvariable}, 
il suffit d'estimer $\nu_{\xi(t)}(V(o,\xi,t))$, ce que nous ferons par la suite.

\begin{lem}\label{ombreclassique}
Il existe une constante $B_0>0$
telle que si $\xi(t)\in\Gamma\widetilde{C}_0$, alors 
$$
\frac{1}{B_0}\,\le\,\nu_{\xi(t)}(V(o,\xi,t))\,\le B_0.
$$
\end{lem}

\noindent
Vu l'encadrement (\ref{simplification1}), 
le lemme \ref{ombreclassique} d\'emontre la partie {\bf a} du th\'eor\`eme, 
avec $A_0=B_0\exp(4\delta \alpha)$.

\begin{demo}Si $\xi(t)\in \Gamma \widetilde{C}_0$, il existe $\gamma\in\Gamma$
tel que $d(\xi(t),\gamma o)\leq \Delta$, avec $\Delta$ le diam\`etre de $C_0$.
Par conformit\'e de $\nu$, ceci implique:
$$
e^{-\delta\Delta}\nu_{\xi(t)}(V(o,\xi,t))\le \nu_{\gamma o}(V(o,\xi,t))=\nu_{o}(V(\gamma^{-1}o,\gamma^{-1}\xi,t))
\le e^{\delta\Delta}\nu_{\xi(t)}(V(o,\xi,t)).
$$
Comme $\nu_o$ est une probabilit\'e, la quantit\'e $\nu_{\xi(t)}(V(o,\xi,t))$ est major\'ee par $\exp(\delta\Delta)$. 
D'autre part, comme la distance de $\gamma^{-1}\xi(t)$ \`a $o$ est inf\'erieure \`a $\Delta$, le lemme \ref{epsilon} donne
$$V(\gamma^{-1}o,\gamma^{-1}\xi,t)=V(\gamma^{-1}\xi(t),\gamma^{-1}\xi,0)\supset B_o(\gamma^{-1}\xi,\varepsilon),$$
et comme $\nu$ est de support $\Lambda_{\Gamma}$, on en d\'eduit:
$$
1\ge \nu_{o}(V(\gamma^{-1}o,\gamma^{-1}\xi,t))\ge \inf_{\eta\in\Lambda_{\Gamma}}\nu_o(B_o(\eta,\varepsilon))=C_{\varepsilon}>0.
$$
\end{demo}

Consid\'erons \`a pr\'esent le cas o\`u $\xi=\xi_i$ est un point parabolique 
de $\Lambda_{\Gamma}$. On notera $(\xi_i(t))_{t\ge 0}$ le rayon $[o\xi_i)$, et
$s_i$ l'instant d'entr\'ee du rayon $[o\xi_i)$ dans l'horoboule ${\cal H}_i$, 
soit encore $\xi_i(s_i)=[o\xi_i)\cap\partial {\cal H}_i$.

\begin{prop}\label{cle}
Il existe des constantes $B_1>0$  et $B_2>0$, telles que pour tout 
 point parabolique $\xi_i$ et tout $t\ge s_i$, on ait:
$$
\frac{1}{B_1} \exp((2\delta_{\Pi}-\delta)(t-s_i))\,\le\,
\nu_{\xi_i(t)}(V(o,\xi_i,t))\,\le \,B_1 \exp((2\delta_{\Pi}-\delta)(t-s_i)),\quad \quad \mbox{et}
$$
$$
\frac{1}{B_2} \exp((2\delta_{\Pi}-\delta)(t-s_i))\,\le\,
\nu_{\xi_i(t)}(\D\setminus V(o,\xi_i,t))\,\le\, B_2  \exp((2\delta_{\Pi}-\delta)(t-s_i)).\quad \quad \quad
$$
\end{prop}

\begin{corol}\label{parabolique}
Cette proposition d\'emontre le th\'eor\`eme quand $\xi=\xi_i\in\Lambda_{pb}$ et $\xi_i(t)\in {\cal H}_i$. 
\end{corol}
\begin{demon}{du corollaire \ref{parabolique}}
En effet, dans ce cas, on a 
$$|d(\xi(t),\Gamma o)-(t-s_i)|=|d(\xi(t),\Gamma o)-d(\xi(t), \Gamma\widetilde{C}_0)|\le \Delta.$$
\end{demon}

\begin{demon}{de la proposition \ref{cle}}
La premi\`ere \'etape est l'estimation des mesures $\nu_{\xi_1(t)}(V(o,\xi_1,t))$ et 
$\nu_{\xi_1(t)}(\D\setminus V(o,\xi_1,t))$ quand $\xi_i=\xi_1$ et $\xi_1(t)\in\widetilde{C}_1$. 
A la fin de la preuve, nous expliquerons comment passer de $\xi_1$ \`a $\xi_i=\gamma_i\xi_1$.

Dans ce premier cas, vus les choix des relev\'es $\widetilde{C}_0$ et $\widetilde{C}_1$ de $C_0$ et $C_1$,
 $\xi_1(s_1)$ est dans l'adh\'erence de $\widetilde{C}_0$, donc \`a distance inf\'erieure 
\`a $\Delta$ de $o$, de sorte que $s_1\le \Delta$. 
On peut donc oublier 
la contribution de $s_1$ dans l'estim\'ee ci-dessus, 
ou plus exactement consid\'erer qu'elle est int\'egr\'ee dans les constantes $B_1$ et $B_2$.

Soit ${\cal D}$ un domaine 
fondamental {\it bor\'elien}  pour l'action de 
$\Pi$ sur $\Lambda_{\Gamma}\setminus\{\xi_1\}$,
i.e. satisfaisant \`a $\nu_x(\cup p{\cal D})= \nu_x(\Lambda_{\Gamma}\setminus\{\xi_1\})$ 
et  $\nu_x({\cal D}\cap p{\cal D})=0$ pour tout $p\neq Id$,. On le choisira 
relativement compact dans $\Lambda_{\Gamma}\setminus\{\xi_1\}$. 
Comme $\nu$ est sans atomes et de support $\Lambda_{\Gamma}$, on a pour tout $x$,
$\nu_x(\Lambda_{\Gamma}\setminus\{\xi_1\})=\nu_x(\D)$, d'o\`u:
$$
\sum_{
\scriptsize
\begin{array}{c}
\scriptsize p \in\Pi,\\
\scriptsize p{\cal D}\subset V(o,\xi_1,t)
\end{array}}\hspace{-0.5cm}\nu_{\xi_1(t)}(p{\cal D})\,
 \leq  \;\nu_{\xi_1(t)}(V(o,\xi_1,t))\;
 \leq \hspace{-0.5cm}
\sum_{
\scriptsize
\begin{array}{c}
\scriptsize p \in\Pi,\\
\scriptsize p{\cal D}\cap V(o,\xi_1,t)\neq\emptyset
\end{array}}\hspace{-0.7cm}\nu_{\xi_1(t)}(p{\cal D}).
$$
En utilisant le lemme \ref{lemme29} avec $K={\cal D}$, on obtient:
\begin{eqnarray}\label{numero1}
\sum_{
\scriptsize
\begin{array}{c}
\scriptsize p \in\Pi,\\
\scriptsize d(o,po)\ge 2t+2K_3
\end{array}}\hspace{-0.7cm}\nu_{\xi_1(t)}(p{\cal D})\,
 \leq  \;\nu_{\xi_1(t)}(V(o,\xi_1,t))\;
 \leq \hspace{-0.5cm}
\sum_{
\scriptsize
\begin{array}{c}
\scriptsize p \in\Pi,\\
\scriptsize d(o,po)\ge 2t-2K_3
\end{array}}\hspace{-0.7cm}\nu_{\xi_1(t)}(p{\cal D}).
\end{eqnarray}

\noindent
De m\^eme, on trouve: 
\begin{eqnarray}\label{numero2}
\sum_{
\scriptsize
\begin{array}{c}
\scriptsize p \in\Pi,\\
\scriptsize d(o,po)\le 2t-2K_3
\end{array}}\hspace{-0.7cm}\nu_{\xi_1(t)}(p{\cal D})
\leq  \;\nu_{\xi_1(t)}(\D\setminus V(o,\xi_1,t))\;
 \leq \hspace{-0.5cm}
\sum_{
\scriptsize
\begin{array}{c}
\scriptsize p \in\Pi,\\
\scriptsize d(o,po)\le 2t+2K_3
\end{array}}\hspace{-0.7cm}\nu_{\xi_1(t)}(p{\cal D})
\end{eqnarray}

\noindent
Aux constantes pr\`es, il nous suffit maintenant d'obtenir 
des estim\'ees des s\'eries 
\begin{eqnarray}\label{deuxseries}
\sum_{
\scriptsize
\begin{array}{c}
\scriptsize p \in\Pi,\\
\scriptsize d(o,po)\ge 2t
\end{array}}\hspace{-0.3cm}\nu_{\xi_1(t)}(p{\cal D})
\quad\quad
\mbox{et}\quad\quad
\sum_{
\scriptsize
\begin{array}{c}
\scriptsize p \in\Pi,\\
\scriptsize d(o,po)\le 2t
\end{array}}\hspace{-0.3cm}\nu_{\xi_1(t)}(p{\cal D})
\end{eqnarray}

\noindent
Nous devons donc calculer  $\nu_{\xi_1(t)}(p{\cal D})$.
Par conformit\'e de $\nu$, nous avons:
$$
\nu_{\xi_1(t)}(p{\cal D})=\int_{\cal D}e^{-\delta\beta_{p\eta}(\xi_1(t),p\xi_1(t))}\,d\nu_{\xi_1(t)}(\eta)\,.
$$
Le lemme \ref{lemme29} fournit une estim\'ee de  $\beta_{p\eta}(\xi_1(t),p\xi_1(t))$ en fonction de $d(o, po)$.
 Si $d(o,po)\le 2t$, nous en d\'eduisons: 
$$
\exp(-2\delta K_3)\le\frac{\nu_{\xi_1(t)}(p{\cal D})}{\nu_{\xi_1(t)}({\cal D})}\le
\exp(2\delta K_3).
$$
Et si $d(o,po)\ge 2t$, alors:
$$
\exp(-2\delta K_3-\delta d(o,po)+2\delta t)\le\frac{\nu_{\xi_1(t)}(p{\cal D})}{\nu_{\xi_1(t)}({\cal D})}\le
\exp(2\delta K_3-\delta d(o,po)+2\delta t).
$$
Dans les deux cas, il reste \`a estimer $\nu_{\xi_1(t)}({\cal D})$. Or par conformit\'e de $\nu$, on a:
\begin{eqnarray}\label{nuxiundet}
\nu_{\xi_1(t)}({\cal D})=\int_{\cal D}e^{-\delta\beta_{\eta}(\xi_1(t),o)}\,d\nu_{o}(\eta)\,.
\end{eqnarray}
Comme au paragraphe \ref{s25}, notons  $y_{\eta}$ l'intersection 
de la g\'eod\'esique $(\eta\xi_1)$ avec $\partial {\cal H}_1$. 
Lorsque $\eta$ varie dans le compact ${\cal D}$ de $\D\setminus\{\xi_1\}$, 
la distance de $o$ \`a la g\'eod\'esique $(\eta\xi_1)$ est inf\'erieure 
\`a la distance de $o$ \`a $y_{\eta}$, 
elle-m\^eme major\'ee par une constante $C=C({\cal D},o,\xi_1)$.
En notant $(\xi_{\eta}(s))_{s\ge 0}$ le rayon $[y_{\eta}\xi_1)$, 
par le lemme \ref{angleausommetbis}, nous avons:
$$
0\le\beta_{\eta}(o,\xi_{\eta}(t))+\beta_{\xi_1}(o,\xi_{\eta}(t))=\beta_{\eta}(o,\xi_{\eta}(t))+t\le 2C.
$$
D'autre part, on a 
$|\beta_{\eta}(\xi_1(t),\xi_{\eta}(t))|\le d(\xi_1(t),\xi_{\eta}(t))\le C$, 
d'o\`u finalement:
$$
t-3C\le \beta_{\eta}(\xi_1(t),o)\le t+C.
$$
En reportant cet encadrement dans le calcul (\ref{nuxiundet}) ci-dessus 
de $\nu_{\xi_1(t)}({\cal D})$, on obtient:
$$
e^{-\delta C-\delta t}\,\nu_o({\cal D})\le \nu_{\xi_1(t)}({\cal D})
\le e^{3\delta C-\delta t}\,\nu_o({\cal D}).
$$

Ainsi, aux constantes multiplicatives pr\`es, 
on est ramen\'e \`a l'estimation des deux s\'eries:
$$
\sum_{d(o,po)\ge 2t}e^{-\delta d(o,po)+\delta t}\quad\mbox{et}\quad
\sum_{d(o,po)\le 2t}\,e^{-\delta t} 
$$
Notons $\displaystyle a_T=\sharp\,\{p\in\Pi,\,d(o,po)\in [T,T+1[\}$. 
Alors $a_T\asymp e^{\delta_{\Pi}T}$ d'apr\`es l'hypoth\`ese (\ref{hypothesecroissance}).
(La notation $f(t)\asymp g(t)$ signifie que pour tout $t\ge 0$, 
le quotient  $f(t)/g(t)$ est compris entre deux constantes strictement positives.)
Alors on a:
$$
\sum_{d(o,po)\ge 2t}e^{-\delta d(o,po)+\delta t}\asymp 
e^{\delta t}\sum_{n=[2t]}^{+\infty}a_n e^{-\delta n}
\asymp e^{\delta t}\sum_{n\ge [2t]}e^{(\delta_{\Pi}-\delta) n}
$$
D'autre part, on sait que $\delta\ge \delta_{\Gamma}$ et $\delta_{\Pi}<\delta_{\Gamma}$.
On en d\'eduit l'estim\'ee voulue:
$$
\nu_{\xi_1(t)}(V(o,\xi_1,t))\asymp  e^{(2\delta_{\Pi}-\delta) t}.
$$
La deuxi\`eme somme est comparable \`a 
$$
e^{-\delta t}\sum_{n=0}^{[2t]}a_n
\asymp e^{-\delta t}\sum_{n=0}^{[2t]}e^{\delta_{\Pi}n}\asymp e^{(2\delta_{\Pi}-\delta) t}, 
$$
d'o\`u 
$$
\nu_{\xi_1(t)}(\D\setminus V(o,\xi_1,t))\asymp  e^{(2\delta_{\Pi}-\delta) t}.
$$

Il reste maintenant \`a consid\'erer le cas o\`u $\xi_i=\gamma_i\xi_1\neq\xi_1$. 
Rappelons la notation $\xi_i(s_i)=[o \,\xi_i(t)[\cap\partial {\cal H}_i$. 
Par choix de $\gamma_i$,  $\gamma_i\widetilde{C}_0$ intersecte le rayon $[o\xi_i)$ 
et son adh\'erence intersecte celle de $\widetilde{C}_i$. 
On en d\'eduit que $\xi_i(s_i)$ est dans le bord  de $\gamma_i\widetilde{C}_0$, 
d'o\`u $d(\xi_i(s_i), \gamma_i o)\le \Delta$.
En utilisant le lemme \ref{encadrement}c (avec $K_2=K_2(\Delta)$), 
et le fait que $V(o,\xi_i,t)=V(\xi_i(s_i),\xi_i,t-s_i)$, on obtient l'encadrement:
\begin{eqnarray}\label{inclusions}
V(\gamma_i o,\xi_i, t-s_i+K_2)\subset
V(o,\xi_i,t)\subset V(\gamma_i o,\xi_i, t-s_i-K_2).
\end{eqnarray}
Notons $(\xi_1(u))_{u\ge 0}$ le rayon $[o\xi_1)$.
Les rayons ${(\xi_i(u+s_i))}_{u\ge 0} $ et ${(\gamma_i\xi_1(u))}_{u\ge 0}$ sont asymptotes, 
d'o\`u pour tout $t\ge s_i$:
$$
d(\xi_i(t),\gamma_i\xi_1(t-s_i))\le d(\xi_i(s_i),\gamma_i\xi_1(0))=d(\xi_i(s_i),\gamma_i o)\le \Delta\,.
$$ 
La conformit\'e de $\nu$ donne pour tout $\eta\in \D$: 
\begin{eqnarray}\label{quot}
\exp(-\delta \Delta)
\le \frac{d\nu_{\xi_i(t)}}{d\nu_{\gamma_i\xi_1(t-s_i)}}(\eta)\le \exp(\delta \Delta).
\end{eqnarray}
Les encadrements (\ref{inclusions}) et (\ref{quot}) montrent 
qu'\`a des constantes multiplicatives pr\`es, il suffit de savoir estimer les quantit\'es  
$$\nu_{\gamma_i\xi_1(t-s_i)}(V(\gamma_i o,\gamma_i\xi_1, t-s_i))=
\nu_{\xi_1(t-s_i)}(V(o,\xi_1, t-s_i)),\quad\mbox{et}
$$ 
$$
\nu_{\gamma_i\xi_1(t-s_i)}(\D\setminus V(\gamma_i o,\gamma_i\xi_1, t-s_i))=
\nu_{\xi_1(t-s_i)}(\D\setminus V(o,\xi_1, t-s_i)).
$$
La premi\`ere partie de la preuve s'applique alors pour donner 
les estim\'ees voulues.
Ceci conclut la d\'emonstration de la proposition. 
\end{demon}

Il reste maintenant \`a traiter le cas g\'en\'eral, 
ce qui est fait dans le lemme ci-dessous et conclut 
la preuve du th\'eor\`eme \ref{courbvariable}: 

\begin{lem}\label{derniercas} 
Il existe une constante $B_3>0$, telle que si $\xi\in \Lambda_{\Gamma}$, 
et $\xi(t)\in \Gamma\widetilde{C}_1$,
alors: 
$$
\frac{1}{B_3}\,\exp\left((2\delta_{\Pi}-\delta)d(\xi(t),\Gamma o)\right)\,\le 
\nu_{\xi(t)}(V(o,\xi,t))\le B_3\,\exp\left((2\delta_{\Pi}-\delta)d(\xi(t),\Gamma o)\right)
$$
\end{lem}

\begin{demo}Soit $\xi\in\Lambda_{\Gamma}$, et $i\in\N$ tel que  $\xi(t)\in\widetilde{C}_i$. 
Nous distinguerons trois cas selon les positions respectives de $\xi$ et $\xi_i$.

\noindent
{\bf Premier cas:} 
$\xi_i\in V(o,\xi,t+K_1+\alpha)$ 
(o\`u $K_1$ est la constante donn\'ee par le lemme \ref{encadrement}b).\\
Montrons que la mesure $\nu_{\xi(t)}(V(o,\xi,t))$ est tr\`es proche de 
$\nu_{\xi_i(t)}(V(o,\xi_i,t))\asymp e^{(2\delta_{\Pi}-\delta)(t-s_i)}$ (proposition \ref{cle}), 
et que la distance de $\xi(t)$ \`a $\Gamma o$ est \`a peu pr\`es $t-s_i$.

\noindent
Le lemme \ref{encadrement} b donne: 
$$
V(o,\xi_i,t+K_1)\subset V(o,\xi,t)\subset V(o,\xi_i,t-K_1).
$$
D'autre part, pour tout $\eta\in\D$, l'in\'egalit\'e triangulaire, 
puis le fait que $\xi_i\in V(o,\xi,t+K_1+\alpha)$ donnent:
$|\beta_{\eta}(\xi(t),\xi_i(t))|\le d(\xi(t),\xi_i(t))\le \alpha$.
Par conformit\'e de $\nu$, on en d\'eduit:
$$
e^{-\delta \alpha}\,\nu_{\xi_i(t)}(V(o,\xi_i,t+K_1))\le \nu_{\xi(t)}(V(o,\xi,t))
\le\,e^{\delta \alpha}\,\nu_{\xi_i(t)}(V(o,\xi_i,t-K_1)).
$$
La premi\`ere estim\'ee de la proposition \ref{cle} donne alors le r\'esultat, puisque
$$
|d(\xi(t),\Gamma o)-(t-s_i)|\le |d(\xi(t),\Gamma o)-d(\xi_i(t),\Gamma o)|+|(t-s_i)-d(\xi_i(t),\Gamma o)|\le \alpha+\Delta.
$$

\noindent
{\bf Deuxi\`eme cas:} 
$\xi_i\notin V(o,\xi,t-K_1-\alpha)$.\\
Le point $\xi_i$ est alors plus proche du point antipodal de 
$\xi$ sur $(o\xi)$, not\'e $\xi'$. Introduisons  alors l'autre intersection $o'$ 
de la 
g\'eod\'esique $(o\xi)$ avec $\partial {\cal H}_i$:
$o'=]\xi(t)\xi)\cap\partial {\cal H}_i$, et la distance $t'=d(\xi(t),o')$. 
On a alors $V(o,\xi,t)=\D\setminus V(o',\xi',t')$, et $\xi_i\in V(o',\xi',t'+K_1+\alpha)$. 
Le lemme \ref{encadrement}b donne 
$$
\D\setminus V(o',\xi_i,t'-K_1)\subset V(o,\xi,t)=\D\setminus V(o',\xi',t')\subset \D\setminus V(o',\xi_i,t'+K_1).
$$
Aux constantes pr\`es, on est donc ramen\'es \`a l'estimation de $\nu_{\xi(t)}(\D\setminus V(o',\xi_i,t'))$.

Comme $o\in \widetilde{C}({\Gamma})$ et $\xi\in\Lambda_{\Gamma}$, 
$o'$ est aussi dans $ \widetilde{C}({\Gamma})$. 
De plus, par d\'efinition, $o'$ appartient au bord de ${\cal H}_i$, 
donc au bord de $\Gamma\widetilde{C}_0$.
On peut donc trouver un $\gamma \in\Gamma$ tel que $d(o',\gamma o)=d(o',\Gamma o)\le \Delta$.
Aux constantes pr\`es, par le lemme \ref{encadrement}c, on est ramen\'e \`a estimer $\nu_{\xi(t)}(\D\setminus V(\gamma o,\xi_i,t'))$.

Notons $(c(s))_{s\ge 0}$ le rayon $[\gamma o\, \xi_i)$. 
En utilisant le fait que $\xi_i\in V(o',\xi',t'+K_1+\alpha)$ et 
que $d(o',\gamma o)\le \Delta$, par les lemmes \ref{encadrement}c et \ref{comparables}, 
on obtient $d(c(t'),\xi(t))\le cste(\alpha,\Delta)$.
Par conformit\'e de $\nu$, on en d\'eduit:
$$
\nu_{\xi(t)}(\D\setminus V(\gamma o,\xi_i,t'))\asymp \nu_{c(t')}(\D\setminus V(\gamma o,\xi_i,t'))
$$
La proposition \ref{cle} donne: 
$\displaystyle 
\nu_{c(t')}(\D\setminus V(\gamma o,\xi_i,t'))\asymp e^{(2\delta_{\Pi}-\delta)t'}.$
En rassemblant les approximations successives ci-dessus, on en d\'eduit 
$$
\nu_{\xi(t)}(V(o,\xi,t))\asymp e^{(2\delta_{\Pi}-\delta)t'}.
$$

Il reste \`a voir que la distance de $\xi(t)$ \`a $\Gamma o$ est \`a peu pr\`es $t'$. 
Ceci d\'ecoule du fait (vu ci-dessus) que $\xi(t)$ est \`a distance born\'ee de $c(t')$, 
et de l'\'egalit\'e $d(c(t'),\Gamma o)=d(c(t'),\gamma o)=t'$.\\

\noindent
{\bf Dernier cas:} $\xi_i\in V(o,\xi,t-K_1-\alpha)\setminus V(o,\xi,t+K_1+\alpha)$. \\
On pose alors $t_1=t-2K_1-2\alpha$ et $t_2=t+2K_1+2\alpha$. On a alors 
$\xi_i\in V(o,\xi,t_1+K_1+\alpha)$ et $\xi_i\notin V(o,\xi,t_2-K_1-\alpha)$.
 D'o\`u le r\'esultat voulu, d'apr\`es 
les deux premiers cas ci-dessus appliqu\'es respectivement \`a $t_1$ et $t_2$.
\end{demo}

Ceci termine donc la preuve du Lemme de l'Ombre en courbure n\'egative variable. 
Nous allons \`a pr\'esent voir comment ces estim\'ees permettent d'\'etudier la non divergence des horosph\`eres.


\section{Moyennes horosph\'eriques}\label{s4}

\subsection{D\'efinitions}\label{41}

Rappelons que $M=\Gamma\backslash \widetilde{M}$ est une vari\'et\'e g\'eom\'etriquement finie \`a courbure major\'ee par $-1$.
Le flot g\'eod\'esique $g=(g^t)_{t\in\R}$ agissant sur $\tm$ 
est un flot hyperbolique, dont les vari\'et\'es fortement instables 
se rel\`event sur $\ttm$ en les vari\'et\'es fortement instables du 
flot g\'eod\'esique $\tilde{g}=(\tilde{g}^t)_{t\in\R}$ de $\ttm$.
Sur $\ttm$, il existe une tr\`es bonne description g\'eom\'etrique de ces ensembles.

Une horosph\`ere $H\subset\M$ centr\'ee en $\xi$ se rel\`eve \`a $\ttm$ en 
une {\it horosph\`ere fortement instable} 
$H^+:=\{u\in\ttm,\,\pi(u)\in H,\,\mbox{et } u^-=\xi \}$. 
C'est encore l'ensemble des vecteurs bas\'es sur $H$ orthogonaux \`a $H$, 
et pointant vers l'ext\'erieur. 
Si $u\in\ttm$, nous noterons respectivement $H(u)\subset M$ l'horosph\`ere de $M$ 
centr\'ee en $u^-$ et passant par le point base $\pi(u)$ de $u$, et $H^+(u)\subset\ttm$ 
l'horosph\`ere fortement instable contenant $u$.
Les horosph\`eres vues sur $\ttm$ sont  les 
vari\'et\'es fortement instables du flot g\'eod\'esique $\tilde{g}$ sur $\ttm$: 
$$
H^+(u)=\widetilde{W}^{su}(u):=\{w\in\ttm,\, \lim_{t\to +\infty}d(\tilde{g}^{-t}u, \tilde{g}^{-t}w)=0\}
$$
De la m\^eme mani\`ere, on d\'efinit l'{\it horosph\`ere fortement stable} $H^-(u)$ de $u$ par:
$$
H^-(u)=\widetilde{W}^{ss}(u):=\{w\in\ttm,\, \lim_{t\to +\infty}d(\tilde{g}^{t}u, \tilde{g}^{t}w)=0\}
$$

\begin{figure}[!ht]
\begin{center}
\input{horosphereND.pstex_t}
\caption{Horosph\`eres et horoboules}
\end{center} 
\end{figure}
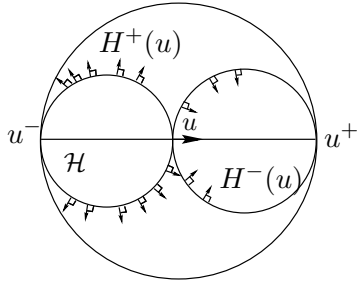


Nous nous servirons de la famille $(d_{H^+})_{H^+\in{\cal H}}$ des {\em distances de Hamenst\"adt} 
sur les horosph\`eres de $\ttm$ 
Si $x\in\M$ est un point quelconque, 
elles sont d\'efinies pour tous $(u,v)\in (H^+)^2$ par
$$
d_{H^+}(u,v)\,=\,\exp\left(\frac{1}{2}\beta_{u^+}(x,u)+\frac{1}{2}\beta_{v^+}(x,v)\right)\,d_x(u^+,v^+)
$$
En fait, nous les consid\'ererons indiff\'eremment comme des distances sur les horosph\`eres de $\M$ et de $\ttm$.
Elles sont bien d\'efinies (car l'expression ci-dessus ne d\'epend pas de $x$), 
elles sont invariantes par isom\'etries au sens o\`u
pour toute isom\'etrie $g$ de $\ttm$, $\displaystyle d_{g H^+}(gu,gv)=d_{H^+}(u,v)$, 
et pouss\'ees par le flot, elles v\'erifient pour tout $t\in\R$:
$$ 
d_{g^t H^+}(g^t u, g^t v)=e^t\,d_{H^+}(u,v).
$$

\begin{figure}[!h]
\begin{center}
\input{distanceND.pstex_t}
\caption{Distance horosph\'erique}\label{distance}
\end{center}
\end{figure}
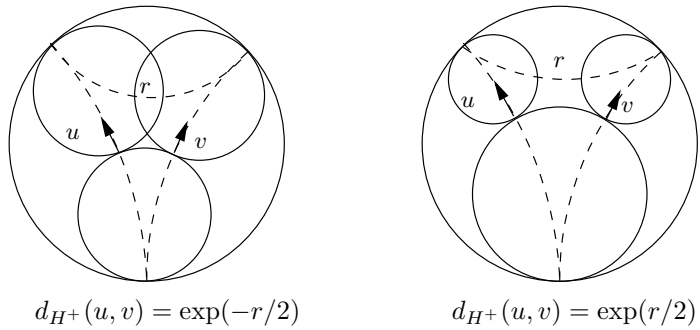

\noindent
G\'eom\'etriquement (figure \ref{distance}), $2\log d_{H^+}(u,v)$ repr\'esente la distance <<~alg\'ebrique~>> 
entre les deux horosph\`eres $H^-(u)$ et $H^-(v)$.

Dans le mod\`ele du demi-plan sup\'erieur de  l'espace hyperbolique $\H$, cette distance s'exprime tr\`es simplement. 
Si $H^+(u)$ est une horosph\`ere horizontale, i.e. si $u^-=+\infty$, et $v$ est un vecteur de $H^+(u)$, 
alors on a $d_{H^+}(u,v)=\frac{d}{h}$, o\`u $d$ est la distance 
euclidienne entre les points base de $u$ et $v$, et $h$ la hauteur euclidienne de l'horosph\`ere $H^+(u)$.\\

%
%

Rappelons que si $\Gamma$ est un groupe g\'eom\'etriquement fini 
v\'erifiant l'hypoth\`ese (\ref{hypothesecroissance}), il est divergent, et que 
par cons\'equent, il existe une  unique densit\'e conforme invariante normalis\'ee $\nu=(\nu_x)_{x\in\M}$
de dimension $\delta_{\Gamma}$ sur $\Lambda_{\Gamma}$ (proposition \ref{uniquedensite}).

Pour pouvoir d\'efinir des moyennes horosph\'eriques au paragraphe suivant, 
nous aurons besoin d'une famille de mesures sur les horosph\`eres fortement instables, 
d\'efinies \`a partir de la mesure de Patterson $\nu_o$.

\begin{prop}\label{mesure}
La famille de mesures d\'efinies sur chaque horosph\`ere $H^+$ par:
$$
d\mu_{H^+}(v)=\exp(\delta_{\Gamma}\beta_{v^+}(o,v))d\nu_o(v^+)
$$ 
est ind\'ependante de $o$, $\Gamma$-invariante au sens o\`u $\gamma_*\mu_{H^+}=\mu_{\gamma H^+}$, 
et passe donc au quotient par $\Gamma$ en une famille de mesures sur 
les vari\'et\'es fortement instables du flot g\'eod\'esique. 
De plus, pouss\'ee par le flot, elle 
v\'erifie pour tout $t\in\R$:
\begin{eqnarray}\label{flot}
dg^{-t}_*\mu_{g^t H^+}(u)=\exp(\delta_{\Gamma}t)\,d\mu_{H^+}(u)
\end{eqnarray}
\end{prop}

Remarquons \'egalement que, la mesure $\nu_o$ ayant pour support $\Lambda_{\Gamma}$, 
chaque mesure $\mu_{H^+}$ d\'efinie ci-dessus a pour support $\{v\in H^+, v^+\in\Lambda_{\Gamma}\}$.
En particulier, si $H^+$ est centr\'ee dans $\Lambda_{\Gamma}$, $\mu_{H^+}$ a pour support $H^+\cap\Omega$.


\subsection{Non divergence des moyennes horosph\'eriques}

Quand $M$ est g\'eom\'etriquement finie, 
la topologie des horosph\`eres fortement instables est bien connue.
Soit  ${\cal E}=\{v\in \tm, v^-\in\Lambda_{\Gamma}\}$. D'apr\`es 
Dal'bo \cite{Dalbo}, 
les feuilles centr\'ees en un point parabolique born\'e sont compactes, 
celles centr\'ees en un point limite radial sont denses dans 
${\cal E}$, et celles qui sont centr\'ees 
hors de l'ensemble limite sont ferm\'ees et non compactes. 
Par cons\'equent, par analogie avec le cas d'un flot, nous appellerons ${\cal E}$ 
l'{\it ensemble non errant du feuilletage horosph\'erique (instable)}. 
Il se d\'ecompose en une union disjointe 
$$
{\cal E}={\cal E}_R\sqcup{\cal E}_{pb}\,,
$$
o\`u ${\cal E}_R$ (resp. ${\cal E}_{pb}$) est l'ensemble des vecteurs de ${\cal E}$ 
centr\'es en un point limite radial (resp. un point parabolique born\'e).


En l'absence de param\'etrisation naturelle des feuilles par un flot,
nous nous 
int\'eressons ici \`a des moyennes sur de grandes boules horosph\'eriques 
pour la mesure $\mu_{H^+}$ d\'efinie au paragraphe pr\'ec\'edent:
 pour toute fonction continue $\psi:\tm\to\R$, et tout $r\ge 0$, nous posons:

$$
M_{r,u}(\psi)=\frac{1}{\mu_{H^+}(B^+(u,r))}\int_{B^+(u,r)}\psi(v)\,d\mu_{H^+}(v).
$$

Si $u\in{\cal E}_{pb}$, les probabilit\'es $(M_{r,u})_{r>0}$ sont \`a support 
compact inclus dans $H^+(u)$. 
Mais si $u\in {\cal E}_R$, la mesure $\mu_{H^+}$ ayant pour support $H^+\cap\Omega$, 
ces mesures $(M_{r,u})_{r>0}$ sont \`a support dans l'ensemble non errant $\Omega$ du flot g\'eod\'esique, 
qui est non compact lorsque $M$ est g\'eom\'etriquement finie avec des cusps.

\noindent
Pour quantifier la non divergence des horosph\`eres, 
nous allons montrer qu'il n'y a pas de perte de masse des mesures $(M_{r,u})_{r>0}$ dans les cusps. 
Plus pr\'ecis\'ement, notre r\'esultat est le suivant:

\begin{theo}\label{nondivergence}Soit $M=\M/\Gamma$ une vari\'et\'e
g\'eom\'etriquement finie telle que tout  sous-groupe parabolique maximal $\Pi$ de $\Gamma$ v\'erifie:
$$
\frac{1}{D}\exp(\delta_{\Pi}T)\le \sharp\{p\in\Pi, d(o,p o)\in [ T,T+1[\}\le D \exp(\delta_{\Pi}T).
$$ 
Fixons $\varepsilon>0$, et $C$ un compact de $\tm$. 
Alors il existe un compact $K_{\varepsilon, C}\subset\Omega$, tel que pour tout $u\in C\cap {\cal E}_R$ et tout $r\ge 0$, 
$$
M_{r,u}(K_{\varepsilon,C})\geq 1-\varepsilon
$$
\end{theo}

\begin{demo}Nous allons \'etudier les moyennes sur le rev\^etement $\ttm$, et 
trouver en fait pour tout $\varepsilon>0$ et $C$ comme ci-dessus un ensemble $\widetilde{K}_{\varepsilon,C}$, 
compact modulo $\Gamma$, qui satisfait l'assertion du th\'eor\`eme ci-dessus. \\

\noindent
{\bf Premi\`ere \'etape:} Dans le cas o\`u $C=\pi^{-1}(C_0)$ est l'ensemble 
des vecteurs bas\'es dans $C_0$, nous allons  introduire le candidat \`a \^etre le compact $K_{\varepsilon, C}$ du th\'eor\`eme, 
puis ramener l'\'enonc\'e ci-dessus \`a la recherche d'estim\'ees de mesures de boules horosph\'eriques.


Comme $\Gamma$ est g\'eom\'etriquement fini, il n'a qu'un nombre fini de cusps, 
on supposera donc comme au paragraphe pr\'ec\'edent qu'il n'en a qu'un seul. 
Rappelons que $\widetilde{C}(\Gamma)$ se d\'ecompose en une union disjointe du cusp et de la partie compacte:
$\widetilde{C}(\Gamma)=\Gamma \widetilde{C}_1\sqcup \Gamma \widetilde{C}_0$, avec
$\widetilde{C}_1={\cal H}_1\cap\widetilde{C}(\Gamma)$, 
${\cal H}_1$ une horoboule centr\'ee au point parabolique born\'e $\xi_1$, 
et  $\Pi$  le stabilisateur de $\xi_1$ dans $\Gamma$.
Nous renvoyons le lecteur aux notations introduites au paragraphe \ref{Lemmedelombre}

\noindent
Notons maintenant ${\cal H}_1^N \subset {\cal H}_1$ l'horoboule << r\'etr\'ecie de $N$ >>, c'est-\`a-dire l'horoboule dont le bord
satisfait $\beta_{\xi_1}(\partial {\cal H}_1,\partial {\cal H}_1^N)=N>0$. 
Le candidat \`a \^etre $\widetilde{K}_{\varepsilon,C_0}$ est 
l'ensemble des vecteurs de $\Lambda_{\Gamma}^2\times\R$ dont le point base appartient au compl\'ementaire 
$\widetilde{C}(\Gamma)\setminus\Gamma {\cal H}_1^N$ de ces horoboules r\'etr\'ecies, pour un $N=N(\varepsilon)$ suffisamment grand.

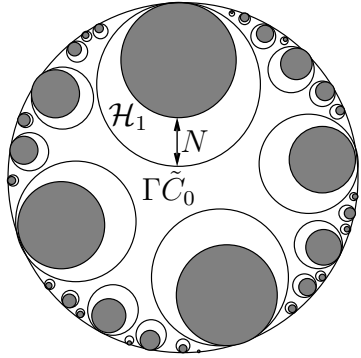
\begin{figure}[!ht]
\begin{center}
\input{horosph_retrecie.pstex_t}
\caption{Horosph\`eres r\'etr\'ecies}
\end{center}
\end{figure}

\noindent
Nous fixons d\'esormais $u\in{\cal E}_R$,  
et nous supposons  que $u$ est bas\'e dans
$\Gamma\widetilde{C}_0$,  tout en gardant 
en m\'emoire que les raisonnements que nous ferons ne doivent d\'ependre que de $C_0$.
Le cas g\'en\'eral d'un compact $C$ quelconque sera trait\'e dans la derni\`ere \'etape.

Nous aurons besoin d'estimer la mesure $\mu_{H^+}(B^+(u,r)\cap{\cal H}_i^N)$, o\`u par abus de notation, $B^+(u,r)\cap{\cal H}_i^N$
d\'esigne l'ensemble des vecteurs de $B^+(u,r)$ dont le point base appartient \`a ${\cal H}_i^N$.
Commen\c cons donc par \'etudier l'intersection g\'eom\'etrique $H^+(u) \cap{\cal H}$, 
o\`u ${\cal H}$ est une horoboule quelconque centr\'ee en un point $\xi\neq u^-$.
Si cette intersection est non vide, 
introduisons le {\it vecteur le plus haut} de $H^+(u)$ dans ${\cal H}$, 
i.e. le vecteur $v$ qui r\'ealise le maximum de $w\to\beta_{\xi}(\partial{\cal H},\pi(w))$ dans  $H^+(u) \cap{\cal H}$.
Par stricte convexit\'e des horosph\`eres, 
ce vecteur est bien d\'efini et s'\'ecrit $v=(u^-,\xi,s(u))$ dans les coordonn\'ees $\ttm\simeq\DDR$.
Nous noterons $h=\beta_{\xi}(\partial {\cal H},\pi(v))$ la hauteur \`a laquelle monte $H^+(u)$ dans ${\cal H}$.

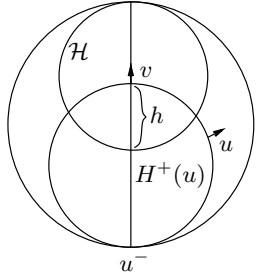
\begin{figure}[!ht]
\begin{center}
\input{lemme43.pstex_t}
\caption{Intersection horosph\`ere-horoboule}
\end{center}
\end{figure}

\begin{lem}\label{boule}Soit ${\cal H}$ une  horoboule centr\'ee en $\xi\in\D$,
et $H^+(u)$ une horosph\`ere fortement instable qui intersecte ${\cal H}$. 
Soient comme ci-dessus $v=(u^-,\xi,s(u))$ le vecteur le plus haut de $H^+(u)$ 
dans ${\cal H}$ et $h=\beta_{\xi}(\partial {\cal H},\pi(v))$.
Si $h\ge \alpha$, on a:
$$
B^+(v,e^{\frac{h-\alpha}{2}})\subset
H^+(u)\cap {\cal H}\subset B^+(v,e^{\frac{h+\alpha}{2}})\,
$$
l'inclusion de droite restant vraie quand $0\le h<\alpha$.
\end{lem}

Dans l'espace hyperbolique $\H$, un calcul donne pr\'ecis\'ement pour tout $h\ge 0$:
$$
 H^+(u)\cap {\cal H}= B^+(v,\sqrt{e^h-1}).
$$

\begin{demon}{du lemme \ref{boule}}
Soit $w\in H^+(u)\cap{\cal H}$. On a alors $0\le \beta_{\xi}(\pi(w),\pi(v))\le h$. 
Si les horosph\`eres (de $\M$) $\pi(H^-(v))$ et $\pi(H^-(w))$  s'intersectent, alors
$$
d_{H^+}(v,w)\le 1\le e^{\frac{h+\alpha}{2}}\,.
$$
On peut donc les supposer disjointes. 
Notons alors ${\cal H}(v)$ (resp. ${\cal H}(w)$) 
l'horoboule dont le bord est $\pi(H^-(v))$ (resp. $\pi(H^-(w))$). 
Dans ce cas, le triangle $(p,q,r)$ int\'erieur \`a $(u^-,\xi,w^+)$ est 
\`a l'ext\'erieur des horoboules ${\cal H}(v)$ et ${\cal H}(w)$.
(En effet, si $r$ appartenait \`a ${\cal H}(v)$, 
par d\'efinition du triangle int\'erieur, 
ceci impliquerait $p\in{\cal H}(v)$, 
puis $q\in{\cal H}(w)$ et enfin $r\in{\cal H}(w)$ 
d'o\`u une contradiction. Les autres cas se traitent de fa\c con analogue.)

Introduisons $a_v$ et $a_w$ les intersections respectives des horosph\`eres $\pi(H^-(v))$ et $\pi(H^-(w))$
avec la g\'eod\'esique $(\xi w^+)=(v^+ w^+)$.
Par d\'efinition de $d_{H^+}$, on a
$\displaystyle d_{H^+}(v,w)=\exp\frac{d(a_v,a_w)}{2}$.
Or
\begin{eqnarray*}
d(a_v,a_w)=\beta_{\xi}(a_w,a_v) &=& \beta_{\xi}(a_w,\pi(w))+\beta_{\xi}(\pi(w),\pi(v))\\
 &\le & d(a_w,\pi(w))+\beta_{\xi}(\pi(w),\pi(v))\\
\end{eqnarray*}
Comme $q$ et $r$ n'appartiennent pas \`a ${\cal H}(w)$, on a $d(a_w,\pi(w))\le d(r,q)\le\alpha$. 
D'autre part, $\pi(w)\in{\cal H}$, d'o\`u
$\beta_{\xi}(\pi(w),\pi(v))\le \beta_{\xi}(\partial{\cal H},\pi(v))=h$. 
Finalement, on en d\'eduit $d(a_v,a_w)=2\log d_{H^+}(v,w)\le h+\alpha$, et 
$$ 
H^+(u)\cap{\cal H}\subset B^+(v,e^{\frac{h+\alpha}{2}})\,.
$$

Supposons maintenant que $\pi(w)\notin{\cal H}$, i.e. $\beta_{\xi}(\pi(w),\pi(v))>h$.
Alors, en reprenant les calculs ci-dessus, on obtient 
$$
d(a_v,a_w)=\beta_{\xi}(a_w,\pi(w))+\beta_{\xi}(\pi(w),\pi(v))>-\alpha+h\,.
$$
Lorsque $-\alpha+h\ge 0$, ceci donne l'autre inclusion 
$\displaystyle 
B^+(v,e^{\frac{h-\alpha}{2}})\subset H^+(u)\cap {\cal H}$.

\begin{figure}[!ht]
\begin{center}
\input{boule_horosph.pstex_t}
\end{center}
\caption{Lemme \ref{boule}}\label{figlemme51}
\end{figure}
\noindent
\end{demon}

Soit maintenant ${\cal H}_i=\gamma_i {\cal H}_1$ une horoboule intersectant $H^+(u)$, 
$v_i=(u^-,\xi_i,s(u))$ le vecteur le plus haut dans l'horoboule, et
$h_i=\beta_{\xi_i}(\partial {\cal H}_i,\pi(v_i))$.
\noindent
On d\'eduit du lemme \ref{boule} que si $h_i\geq N+\alpha$, alors 
$$
B^+(v_i,e^{\frac{h_i-N-\alpha}{2}})\subset H^+(u)\cap {\cal H}_i^N
\subset B^+(v_i,e^{\frac{h_i-N+\alpha}{2}})\,
$$
l'inclusion de droite restant vraie pour $h_i\ge N$.
Si $h_i\leq N$, on a trivialement 
 $H^+(u)\cap {\cal H}_i^N=\emptyset$.

On ne regarde donc que les $h_i\geq N$.
Rappelons que nous cherchons $N=N(\varepsilon)>0$ assez grand pour que, pour tout $u\in{\cal E}_R\cap C_0$, 
\begin{eqnarray}\label{amajorer}
\frac{\mu_{H^+}(B^+(u,r)\cap\Gamma{\cal H}_1^N)}{\mu_{H^+}(B^+(u,r))}\le\varepsilon\,.
\end{eqnarray}
Introduisons l'ensemble $I_{u,r,N}$ des $i\in\N$, 
tels que $B^+(u,r)\cap{\cal H}_i^N\neq\emptyset$. 
(On a donc $h_i\ge N$ pour $i\in I_{u,r,N}$.)
A l'aide du lemme \ref{boule}, on peut majorer le num\'erateur de (\ref{amajorer}) ci-dessus par:
$$
\sum_{i\in I_{u,r,N}} \mu_{H^+}(B^+(v_i,e^{\frac{h_i-N+\alpha}{2}}))\,.
$$
Quitte \`a remplacer $N$ par $N-\alpha$, nous oublierons d\'esormais la constante $\alpha$ dans l'expression ci-dessus.

Pour minorer le d\'enominateur, 
rappelons que les horoboules ${\cal H}_i$ sont deux \`a deux disjointes, et donc les 
boules $B^+(v_i,e^{h_i/2})$ aussi.
Ce d\'enominateur est donc minor\'e par 
$$
\sum_{i\in I_{u,r,N}} \mu_{H^+}(B^+(u,r)\cap B^+(v_i,e^{h_i/2}))\,.
$$
Dans la majoration du num\'erateur, les boules $B^+(u,r)$ n'interviennent plus que via l'ensemble d'indices $I_{u,r,N}$.
Pour le d\'enominateur, le probl\`eme est qu'il n'y a aucune 
raison pour que les boules $B^+(v_i,e^{h_i/2})$ 
soient totalement incluses dans $B^+(u,r)$.
En revanche, 
 comme $u\in\Gamma\widetilde{C}_0$, $u$ n'appartient \`a aucune horoboule ${\cal H}_i$, 
et si $i\in I_{u,r,N}$, alors $B^+(u,r)$ intersecte $B^+(v_i, e^{\frac{h_i-N}{2}})$, d'o\`u:
$$
e^{h_i/2}\le d_{H^+}(u,v_i)\le r+e^{\frac{h_i-N}{2}}\,.
$$
Si $N\ge 2\log 3$, on en d\'eduit pour tout $i\in I_{u,r,N}$:
$$
B^+(v_i,e^{h_i/2})\subset B^+(u,3r).
$$
Le d\'enominateur est donc minor\'e par:
$$
\frac{\mu_{H^+}(B^+(u,r))}{\mu_{H^+}(B^+(u,3r))}\,\times\,\sum_{i\in I_{u,r,N}} \mu_{H^+}(B^+(v_i,e^{\frac{h_i-N}{2}}))\,.
$$
Pour finir, on obtient la majoration suivante (pour $N\ge 2\log 3$):
$$
\frac{\mu_{H^+}(B^+(u,r)\cap\Gamma{\cal H}_1^N)}{\mu_{H^+}(B^+(u,r))}\le
\frac{\sum_{i\in I_{u,r,N}} \mu_{H^+}(B^+(v_i,e^{\frac{h_i-N}{2}}))}{\sum_{i\in I_{u,r,N}} \mu_{H^+}(B^+(v_i,e^{h_i/2}))}
\,\times\,\frac{\mu_{H^+}(B^+(u,3r))}{\mu_{H^+}(B^+(u,r))}.
$$

Dans la deuxi\`eme  \'etape, nous majorons le premier quotient par $\varepsilon$ pour $N$ assez grand. 
A la troisi\`eme \'etape, nous montrons que le deuxi\`eme quotient 
est born\'e uniform\'ement en $u\in{\cal E}_R\cap \Gamma\widetilde{C}_0$ et $r>0$.
Pour finir, la quatri\`eme \'etape traitera le cas d'un compact $C$ quelconque diff\'erent de $C_0$.\\

\noindent
{\bf Deuxi\`eme \'etape:} Pour tout $\varepsilon>0$, il existe $N=N(\varepsilon)>0$, 
tel que pour tout $u\in{\cal E}_R\cap \Gamma\widetilde{C}_0$, tout $r>0$ et $i\in I_{u,r,N}$, on a:
\begin{eqnarray}\label{deuxiemeetape}
\frac{\mu_{H^+}(B^+(v_i,e^{(h_i-N)/2}))}{\mu_{H^+}(B^+(v_i,e^{h_i/2}))}\,\le\,\varepsilon.
\end{eqnarray}

\noindent
Notons $w_i=g^{-h_i}v_i\in\partial {\cal H}_i$. D'apr\`es la relation (\ref{flot}) du paragraphe \ref{41},
la quantit\'e ci-dessus est encore \'egale \`a:

$$
\frac{\mu_{H^+ (w_i)}(B^+(w_i,e^{(-h_i-N)/2}))}{\mu_{H^+ (w_i)}( B^+(w_i,e^{-h_i/2}))}.
$$

\noindent
Remarquons qu'\`a $i$ fix\'e, la quantit\'e 
ci-dessus tend clairement vers $0$ quand $N\to\infty$; 
mais la difficult\'e vient du fait qu'on cherche une uniformit\'e en $i\in I_{u,r,N}$.

\noindent
Pour tout $i$, notons $x_i\in\M$ le point base de $w_i$, 
et $\xi_i=\gamma_i\xi_1\in\Lambda_{pb}$ le point du bord vers lequel il pointe. 
Remarquons que $x_i\in\widetilde{C}(\Gamma)$,
puisque $w_i\in\Lambda_{\Gamma}^2\times\R$.
Rappelons que la mesure  $\mu_{H^+}$ est d\'efinie pour $x\in\M$ quelconque par:
$$
d\mu_{H^+}(v)=\exp(\delta_{\Gamma}\beta_{v^+}(x,v))\,d\nu_x(v^+),
$$
avec $\nu_x$ la mesure de Patterson sur $\D$ vue du point $x$. 
En particulier, si $x=x_i$, la quantit\'e $\exp(\delta_{\Gamma}\beta_{v^+}(x_i,v))$ 
est born\'ee par $\exp\delta_{\Gamma}d(x_i,\pi(v))\le e^{\delta_{\Gamma}\alpha}$  
pour tout $v\in B^+(w_i,e^{-h_i/2})\subset B^+(w_i,1)$.

\noindent
La quantit\'e \`a estimer est donc uniform\'ement proche de 
$$
\frac{\nu_{x_i}(B^+(w_i,e^{(-h_i-N)/2}))}{\nu_{x_i}( B^+(w_i,e^{-h_i/2}))}\,,
$$
o\`u par abus de notation, $B^+(w_i,e^{(\dots)/2})$ d\'esigne l'image 
dans $\D$ de la boule horosph\'erique par la projection naturelle de 
$H^+(w_i)$ dans $\D$ qui \`a $w=(u^-,w^+,s(w_i))$ associe $w^+$.

\noindent
Soient $x\in\M, \xi\in\D$ et $t\in\R^+$; rappelons que
$V(x,\xi,t)\subset\D$ est l'ensemble des points $\eta\in\D$ dont 
le projet\'e sur le rayon g\'eod\'esique $[x,\xi)$ est \`a distance sup\'erieure \`a $t$ de $x$. 
Nous avons alors:

\begin{lem}\label{projection}Si on note $P_w$ la projection  naturelle  de $H^+(w)$ dans $\D\setminus\{w^-\}$,
et $x=\pi(w)$ le point base de $w$ dans $\M$,
on a pour tout $s\geq \alpha$:
$$
V(x,w^+,s+\alpha)\subset P_w(B^+(w,e^{-s}))\subset V(x,w^+,s-\alpha).
$$
\end{lem}

\begin{demo}Ecrivons $w=(w^-,w^+,s(w))$.  Soit $s\ge 0$ et $z\in B^+(w,e^{-s})$.
Consid\'erons le triangle id\'eal $(w^-,w^+,z^+)$, et son triangle int\'erieur $(p,q,r)$, 
dont le diam\`etre est born\'e par $\alpha$. 

\begin{figure}[!ht]
\begin{center}
\input{lemme42.pstex_t}
\caption{Lemme \ref{projection}}
\end{center}
\end{figure}

\noindent
Notons $a_w=(w^+z^+)\cap H^-(z)$ et $a_z= (w^+z^+)\cap H^-(w)$. 
On a $d_{H^+}(w,z)=e^{-d(a_w,a_z)/2}\le e^{-s}$, d'o\`u $d(a_w,a_z)\ge 2s$. 
D'autre part, en observant les appartenances des diff\'erents points aux diverses horosph\`eres, on trouve 
\begin{eqnarray*}
d(a_w,r)&=& \beta_{z^+}(a_w,r)=\beta_{z^+}(\pi(z),q)=d(\pi(z),q)\\
&=&\beta_{w^-}(q,\pi(z))=\beta_{w^-}(p,\pi(w))=d(w,p)\\
&=&\beta_{w^+}(\pi(w),p)=\beta_{w^+}(a_z,r)=d(a_z,r)\,
\end{eqnarray*}
donc toutes ces distances sont \'egales \`a $d(a_w,a_z)/2\ge s$. 
Le triangle $(w^-,w^+,z^+)$ \'etant fin (proposition \ref{trianglesfins}), 
le projet\'e de $z^+$ sur $[w w^+)$ est \`a distance inf\'erieure ou \'egale \`a $\alpha$ de $p$,
et $d(w,p)\ge s$ d'apr\`es ci-dessus. 
Donc $z\in V(x,w^+,s-\alpha)$.
Le m\^eme type d'arguments montre que si $z\in V(x,w^+,s+\alpha)$, alors $d_{H^+}(w,z)\le \exp(-s)$.
\end{demo}

\noindent
Nous sommes donc ramen\'es \`a estimer la quantit\'e:
$$
\frac{\nu_{x_i}(V(x_i,\xi_i,\frac{h_i+N}{2}))}{\nu_{x_i}(V(x_i,\xi_i,\frac{h_i}{2}))}.
$$
Le point base $x_i$ de $w_i$ appartient \`a $\partial {\cal H}_i\cap\widetilde{C}(\Gamma)$. 
Il est donc dans le bord de $\gamma\widetilde{C}_0$, pour un certain $\gamma\in\Gamma$ tel 
que $\gamma{\cal H}_1={\cal H}_i$, et donc $\gamma\xi_1=\xi_i$. On en d\'eduit $d(x_i,\gamma o)\le \mbox{diam}(\widetilde{C}_0)=\Delta$.
En utilisant la $\Gamma$-invariance et la conformit\'e de $\nu=(\nu_x)_{x\in\M}$, 
ainsi que le lemme \ref{encadrement}c, on trouve que la quantit\'e ci-dessus est proche, \`a des constantes uniformes pr\`es, de:
$$
\frac{\nu_{o}(V(o,\xi_1,\frac{h_i+N}{2}))}{\nu_{o}(V(o,\xi_1,\frac{h_i}{2}))} .
$$

\noindent
A l'aide du th\'eor\`eme \ref{courbvariable}, on peut estimer cette quantit\'e:
$$
\frac{\nu_{o}(V(o,\xi_1,\frac{h_i+N}{2}))}{\nu_{o}(V(o,\xi_1,\frac{h_i}{2}))}\le
A_1^2\exp\left(2(\delta_{\Pi}-\delta_{\Gamma})(\frac{h_i+N}{2}-\frac{h_i}{2} )\right)=
\,A_1^2\exp (\delta_{\Pi}-\delta_{\Gamma})N.
$$
Ceci tend vers $0$ uniform\'ement en $i$ quand $N\to +\infty$. 
(Remarquons qu'on n'a utilis\'e jusqu'\`a maintenant qu'une partie 
du th\'eor\`eme \ref{courbvariable}, la proposition \ref{cle}.)\\

\noindent
{\bf Troisi\`eme \'etape: }
\begin{lem}\label{doublante}
Il existe une constante $C>0$, telle que pour tout $u\in {\cal E}_R$ et tout $r>0$, on a:
$$
\frac{\mu_{H^+}(B^+(u,3r))}{\mu_{H^+}(B^+(u,r))}\leq C\,.
$$
\end{lem}

\begin{demo}Comme $u^-\in\Lambda_R$, on peut d\'efinir $s\ge 0$ le plus petit 
r\'eel sup\'erieur \`a $\log 3r$, tel que $g^{-s}u\in\Gamma\widetilde{C}_0$. 
On a alors:
$$
\frac{\mu_{H^+}(B^+(u,3r))}{\mu_{H^+}(B^+(u,r))}=
\frac{\mu_{H^+}(B^+(g^{-s}u,3re^{-s}))}{\mu_{H^+ }(B^+(g^{-s}u,re^{-s}))}.
$$
Or il existe $\gamma\in\Gamma$, tel que $d(\gamma o,g^{-s }u)\le \Delta$. 
Par les lemmes \ref{projection}, \ref{encadrement}b, la conformit\'e  et l'invariance de $\nu$ par $\Gamma$, 
la quantit\'e ci-dessus est donc, \`a des constantes uniformes pr\`es, proche de:
$$
\frac{\nu_o(V(o,\gamma^{-1}u^+, s-\log r-\log 3))}{\nu_o(V(o,\gamma^{-1}u^+, s-\log r))} .
$$
Les deux points du rayon $[o\,\gamma^{-1}u^+)$ \`a distance respective $s-\log r-\log 3$ et $s-\log r$ de $o$ 
\'etant \`a distance $\log 3$ l'un de l'autre, 
on peut supposer qu'ils sont simultan\'ement soit dans $\Gamma\widetilde{C}_0$, soit dans $\Gamma\widetilde{C}_1$.
Dans le premier cas, d'apr\`es le th\'eor\`eme \ref{courbvariable}, 
la quantit\'e ci-dessus est major\'ee par $\displaystyle A_0^2\exp(\delta_{\Gamma }\log 3 )$. 
Dans le second cas, elle est inf\'erieure \`a 
$$
A_1^2\exp(\delta_{\Gamma}\log 3+|2\delta_{\Pi}-\delta_{\Gamma}|\log 3).
$$ Ceci conclut la preuve du lemme.
\end{demo}

\noindent
{\bf Quatri\`eme et derni\`ere \'etape:} Il reste pour d\'emontrer le 
th\'eor\`eme \`a traiter le cas d'un compact $C$ quelconque de $\Omega$. 
Quitte \`a le d\'ecomposer en une union de plus petits compacts, on peut 
supposer qu'il existe un $T_0>0$, tel que $\pi(g^{-T_0}C)\subset \widetilde{C}_0$.
Alors, il est facile de voir que le compact $K_{\varepsilon,C}=g^{T_0}(K_{\varepsilon,C_0})$ convient.
\end{demo}



\end{document}

%% file: coordo_bordND.pstex_t
\begin{picture}(0,0)%
\epsfig{file=coordo_bordND.pstex}%
\end{picture}%
\setlength{\unitlength}{3947sp}%
\begingroup\makeatletter\ifx\SetFigFont\undefined%
\gdef\SetFigFont#1#2#3#4#5{%
  \reset@font\fontsize{#1}{#2pt}%
  \fontfamily{#3}\fontseries{#4}\fontshape{#5}%
  \selectfont}%
\fi\endgroup%
\begin{picture}(1729,1638)(1,-802)
\put(1173,-311){\makebox(0,0)[lb]{\smash{\SetFigFont{10}{12.0}{\rmdefault}{\mddefault}{\updefault}$o$}}}
\put(  1,-26){\makebox(0,0)[lb]{\smash{\SetFigFont{10}{12.0}{\rmdefault}{\mddefault}{\updefault}$v^-$}}}
\put(786,  6){\makebox(0,0)[lb]{\smash{\SetFigFont{10}{12.0}{\rmdefault}{\mddefault}{\updefault}$v$}}}
\put(787,764){\makebox(0,0)[lb]{\smash{\SetFigFont{7}{8.4}{\rmdefault}{\mddefault}{\updefault}$\partial\widetilde{M}$}}}
\put(402,-543){\makebox(0,0)[lb]{\smash{\SetFigFont{10}{12.0}{\rmdefault}{\mddefault}{\updefault}$|\beta_{v^-}(\pi(v),o)|$}}}
\put(1730,-88){\makebox(0,0)[lb]{\smash{\SetFigFont{10}{12.0}{\rmdefault}{\mddefault}{\updefault}$v^+$}}}
\end{picture}

%% file: triangle_int.pstex_t
\begin{picture}(0,0)%
\epsfig{file=triangle_int.pstex}%
\end{picture}%
\setlength{\unitlength}{3947sp}%
\begingroup\makeatletter\ifx\SetFigFont\undefined%
\gdef\SetFigFont#1#2#3#4#5{%
  \reset@font\fontsize{#1}{#2pt}%
  \fontfamily{#3}\fontseries{#4}\fontshape{#5}%
  \selectfont}%
\fi\endgroup%
\begin{picture}(2274,1949)(1051,-1411)
\put(1051,-1411){\makebox(0,0)[lb]{\smash{\SetFigFont{12}{14.4}{\rmdefault}{\mddefault}{\updefault}$b$}}}
\put(2294,389){\makebox(0,0)[lb]{\smash{\SetFigFont{12}{14.4}{\rmdefault}{\mddefault}{\updefault}$a$}}}
\put(3322,-918){\makebox(0,0)[lb]{\smash{\SetFigFont{12}{14.4}{\rmdefault}{\mddefault}{\updefault}$c$}}}
\put(2058,-339){\makebox(0,0)[lb]{\smash{\SetFigFont{12}{14.4}{\rmdefault}{\mddefault}{\updefault}$r$}}}
\put(2401,-768){\makebox(0,0)[lb]{\smash{\SetFigFont{12}{14.4}{\rmdefault}{\mddefault}{\updefault}$p$}}}
\put(2572,-339){\makebox(0,0)[lb]{\smash{\SetFigFont{12}{14.4}{\rmdefault}{\mddefault}{\updefault}$q$}}}
\end{picture}

%% file: voisinage.pstex_t
\begin{picture}(0,0)%
\epsfig{file=voisinage.pstex}%
\end{picture}%
\setlength{\unitlength}{3947sp}%
\begingroup\makeatletter\ifx\SetFigFont\undefined%
\gdef\SetFigFont#1#2#3#4#5{%
  \reset@font\fontsize{#1}{#2pt}%
  \fontfamily{#3}\fontseries{#4}\fontshape{#5}%
  \selectfont}%
\fi\endgroup%
\begin{picture}(2767,1508)(91,-808)
\put(2261,-362){\makebox(0,0)[lb]{\smash{\SetFigFont{11}{13.2}{\rmdefault}{\mddefault}{\updefault}$\eta$}}}
\put(2282,-50){\makebox(0,0)[lb]{\smash{\SetFigFont{11}{13.2}{\rmdefault}{\mddefault}{\updefault}$\xi$}}}
\put( 91,-68){\makebox(0,0)[lb]{\smash{\SetFigFont{11}{13.2}{\rmdefault}{\mddefault}{\updefault}$x$}}}
\put(1053, 37){\makebox(0,0)[lb]{\smash{\SetFigFont{11}{13.2}{\rmdefault}{\mddefault}{\updefault}$\xi_x(t)$}}}
\put(2858,-102){\makebox(0,0)[lb]{\smash{\SetFigFont{11}{13.2}{\rmdefault}{\mddefault}{\updefault}$V(x,\xi,t)$}}}
\end{picture}

%% file: lemme29.pstex_t
\begin{picture}(0,0)%
\epsfig{file=lemme29.pstex}%
\end{picture}%
\setlength{\unitlength}{3947sp}%
\begingroup\makeatletter\ifx\SetFigFont\undefined%
\gdef\SetFigFont#1#2#3#4#5{%
  \reset@font\fontsize{#1}{#2pt}%
  \fontfamily{#3}\fontseries{#4}\fontshape{#5}%
  \selectfont}%
\fi\endgroup%
\begin{picture}(2740,2677)(226,-1817)
\put(1665,725){\makebox(0,0)[lb]{\smash{\SetFigFont{9}{10.8}{\rmdefault}{\mddefault}{\updefault}$\xi$}}}
\put(548,-795){\makebox(0,0)[lb]{\smash{\SetFigFont{9}{10.8}{\rmdefault}{\mddefault}{\updefault}$\eta$}}}
\put(226,-795){\makebox(0,0)[lb]{\smash{\SetFigFont{9}{10.8}{\rmdefault}{\mddefault}{\updefault}$K$}}}
\put(1493,-404){\makebox(0,0)[lb]{\smash{\SetFigFont{9}{10.8}{\rmdefault}{\mddefault}{\updefault}$\xi(t)$}}}
\put(2966,-865){\makebox(0,0)[lb]{\smash{\SetFigFont{9}{10.8}{\rmdefault}{\mddefault}{\updefault}$p\eta$}}}
\put(2190,-668){\makebox(0,0)[lb]{\smash{\SetFigFont{9}{10.8}{\rmdefault}{\mddefault}{\updefault}$py_{\eta}$}}}
\put(2646,-348){\makebox(0,0)[lb]{\smash{\SetFigFont{9}{10.8}{\rmdefault}{\mddefault}{\updefault}$po$}}}
\put(2167,-51){\makebox(0,0)[lb]{\smash{\SetFigFont{9}{10.8}{\rmdefault}{\mddefault}{\updefault}$p\xi(t)$}}}
\put(942,-616){\makebox(0,0)[lb]{\smash{\SetFigFont{9}{10.8}{\rmdefault}{\mddefault}{\updefault}$y_{\eta}$}}}
\put(1085,-1069){\makebox(0,0)[lb]{\smash{\SetFigFont{9}{10.8}{\rmdefault}{\mddefault}{\updefault}$o$}}}
\end{picture}

%% file: horosphereND.pstex_t
\begin{picture}(0,0)%
\epsfig{file=horosphereND.pstex}%
\end{picture}%
\setlength{\unitlength}{3947sp}%
\begingroup\makeatletter\ifx\SetFigFont\undefined%
\gdef\SetFigFont#1#2#3#4#5{%
  \reset@font\fontsize{#1}{#2pt}%
  \fontfamily{#3}\fontseries{#4}\fontshape{#5}%
  \selectfont}%
\fi\endgroup%
\begin{picture}(1992,1750)(601,-921)
\put(601,-51){\makebox(0,0)[lb]{\smash{\SetFigFont{11}{13.2}{\rmdefault}{\mddefault}{\updefault}$u^-$}}}
\put(2593,-75){\makebox(0,0)[lb]{\smash{\SetFigFont{11}{13.2}{\rmdefault}{\mddefault}{\updefault}$u^+$}}}
\put(1934,-332){\makebox(0,0)[lb]{\smash{\SetFigFont{11}{13.2}{\rmdefault}{\mddefault}{\updefault}$H^-(u)$}}}
\put(1177,521){\makebox(0,0)[lb]{\smash{\SetFigFont{11}{13.2}{\rmdefault}{\mddefault}{\updefault}$H^+(u)$}}}
\put(956,-227){\makebox(0,0)[lb]{\smash{\SetFigFont{10}{12.0}{\rmdefault}{\mddefault}{\updefault}${\cal H}$}}}
\put(1701, 28){\makebox(0,0)[lb]{\smash{\SetFigFont{10}{12.0}{\rmdefault}{\mddefault}{\updefault}$u$}}}
\end{picture}

%% file: distanceND.pstex_t
\begin{picture}(0,0)%
\epsfig{file=distanceND.pstex}%
\end{picture}%
\setlength{\unitlength}{3947sp}%
\begingroup\makeatletter\ifx\SetFigFont\undefined%
\gdef\SetFigFont#1#2#3#4#5{%
  \reset@font\fontsize{#1}{#2pt}%
  \fontfamily{#3}\fontseries{#4}\fontshape{#5}%
  \selectfont}%
\fi\endgroup%
\begin{picture}(4338,2015)(100,-1276)
\put(2960, 99){\makebox(0,0)[lb]{\smash{\SetFigFont{9}{10.8}{\rmdefault}{\mddefault}{\updefault}$u$}}}
\put(3969, 70){\makebox(0,0)[lb]{\smash{\SetFigFont{9}{10.8}{\rmdefault}{\mddefault}{\updefault}$v$}}}
\put(480,-102){\makebox(0,0)[lb]{\smash{\SetFigFont{9}{10.8}{\rmdefault}{\mddefault}{\updefault}$u$}}}
\put(1289,-160){\makebox(0,0)[lb]{\smash{\SetFigFont{9}{10.8}{\rmdefault}{\mddefault}{\updefault}$v$}}}
\put(942,185){\makebox(0,0)[lb]{\smash{\SetFigFont{9}{10.8}{\rmdefault}{\mddefault}{\updefault}$r$}}}
\put(284,-1231){\makebox(0,0)[lb]{\smash{\SetFigFont{10}{12.0}{\rmdefault}{\mddefault}{\updefault}$d_{H^+}(u,v)=\exp(-r/2)$}}}
\put(2905,-1231){\makebox(0,0)[lb]{\smash{\SetFigFont{10}{12.0}{\rmdefault}{\mddefault}{\updefault}$d_{H^+}(u,v)=\exp(r/2)$}}}
\put(3542,352){\makebox(0,0)[lb]{\smash{\SetFigFont{9}{10.8}{\rmdefault}{\mddefault}{\updefault}$r$}}}
\end{picture}

%% file: horosph_retrecie.pstex_t
\begin{picture}(0,0)%
\epsfig{file=horosph_retrecie.pstex}%
\end{picture}%
\setlength{\unitlength}{3947sp}%
\begingroup\makeatletter\ifx\SetFigFont\undefined%
\gdef\SetFigFont#1#2#3#4#5{%
  \reset@font\fontsize{#1}{#2pt}%
  \fontfamily{#3}\fontseries{#4}\fontshape{#5}%
  \selectfont}%
\fi\endgroup%
\begin{picture}(2219,2216)(159,-1468)
\put(1019,-491){\makebox(0,0)[lb]{\smash{\SetFigFont{12}{14.4}{\rmdefault}{\mddefault}{\updefault}$\Gamma\tilde{C}_0$}}}
\put(1262,-188){\makebox(0,0)[lb]{\smash{\SetFigFont{12}{14.4}{\rmdefault}{\mddefault}{\updefault}$N$}}}
\put(811,-19){\makebox(0,0)[lb]{\smash{\SetFigFont{12}{14.4}{\rmdefault}{\mddefault}{\updefault}${\cal H}_1$}}}
\end{picture}

%% file: lemme43.pstex_t
\begin{picture}(0,0)%
\epsfig{file=lemme43.pstex}%
\end{picture}%
\setlength{\unitlength}{3947sp}%
\begingroup\makeatletter\ifx\SetFigFont\undefined%
\gdef\SetFigFont#1#2#3#4#5{%
  \reset@font\fontsize{#1}{#2pt}%
  \fontfamily{#3}\fontseries{#4}\fontshape{#5}%
  \selectfont}%
\fi\endgroup%
\begin{picture}(1560,1699)(9,-939)
\put(724,-939){\makebox(0,0)[lb]{\smash{\SetFigFont{9}{10.8}{\rmdefault}{\mddefault}{\updefault}$u^-$}}}
\put(917, -7){\makebox(0,0)[lb]{\smash{\SetFigFont{9}{10.8}{\rmdefault}{\mddefault}{\updefault}$h$}}}
\put(1351,-184){\makebox(0,0)[lb]{\smash{\SetFigFont{9}{10.8}{\rmdefault}{\mddefault}{\updefault}$u$}}}
\put(403,394){\makebox(0,0)[lb]{\smash{\SetFigFont{9}{10.8}{\rmdefault}{\mddefault}{\updefault}${\cal H}$}}}
\put(853,282){\makebox(0,0)[lb]{\smash{\SetFigFont{9}{10.8}{\rmdefault}{\mddefault}{\updefault}$v$}}}
\put(821,-377){\makebox(0,0)[lb]{\smash{\SetFigFont{9}{10.8}{\rmdefault}{\mddefault}{\updefault}$H^+(u)$}}}
\end{picture}

%% file: boule_horosph.pstex_t
\begin{picture}(0,0)%
\epsfig{file=boule_horosph.pstex}%
\end{picture}%
\setlength{\unitlength}{3947sp}%
\begingroup\makeatletter\ifx\SetFigFont\undefined%
\gdef\SetFigFont#1#2#3#4#5{%
  \reset@font\fontsize{#1}{#2pt}%
  \fontfamily{#3}\fontseries{#4}\fontshape{#5}%
  \selectfont}%
\fi\endgroup%
\begin{picture}(3223,3104)(76,-2357)
\put(1740,-2357){\makebox(0,0)[lb]{\smash{\SetFigFont{10}{12.0}{\rmdefault}{\mddefault}{\updefault}$u^-$}}}
\put(3141,  6){\makebox(0,0)[lb]{\smash{\SetFigFont{10}{12.0}{\rmdefault}{\mddefault}{\updefault}$w^+$}}}
\put(1275,152){\makebox(0,0)[lb]{\smash{\SetFigFont{10}{12.0}{\rmdefault}{\mddefault}{\updefault}$H^-(v)$}}}
\put(953,-606){\makebox(0,0)[lb]{\smash{\SetFigFont{10}{12.0}{\rmdefault}{\mddefault}{\updefault}$v$}}}
\put(2324,-432){\makebox(0,0)[lb]{\smash{\SetFigFont{10}{12.0}{\rmdefault}{\mddefault}{\updefault}$w$}}}
\put(2266,-986){\makebox(0,0)[lb]{\smash{\SetFigFont{10}{12.0}{\rmdefault}{\mddefault}{\updefault}$q$}}}
\put(1275,-957){\makebox(0,0)[lb]{\smash{\SetFigFont{10}{12.0}{\rmdefault}{\mddefault}{\updefault}$p$}}}
\put(574,-1016){\makebox(0,0)[lb]{\smash{\SetFigFont{10}{12.0}{\rmdefault}{\mddefault}{\updefault}${\cal H}$}}}
\put(1808,-213){\makebox(0,0)[lb]{\smash{\SetFigFont{10}{12.0}{\rmdefault}{\mddefault}{\updefault}$r$}}}
\put(1462,-181){\makebox(0,0)[lb]{\smash{\SetFigFont{10}{12.0}{\rmdefault}{\mddefault}{\updefault}$a_v$}}}
\put(2029,-181){\makebox(0,0)[lb]{\smash{\SetFigFont{10}{12.0}{\rmdefault}{\mddefault}{\updefault}$a_w$}}}
\put(2753,-1221){\makebox(0,0)[lb]{\smash{\SetFigFont{10}{12.0}{\rmdefault}{\mddefault}{\updefault}$H^+(u)$}}}
\put( 76,-55){\makebox(0,0)[lb]{\smash{\SetFigFont{10}{12.0}{\rmdefault}{\mddefault}{\updefault}$\xi=v^+$}}}
\put(2785,-811){\makebox(0,0)[lb]{\smash{\SetFigFont{10}{12.0}{\rmdefault}{\mddefault}{\updefault}$H^-(w)$}}}
\end{picture}

%% file: lemme42.pstex_t
\begin{picture}(0,0)%
\epsfig{file=lemme42.pstex}%
\end{picture}%
\setlength{\unitlength}{3947sp}%
\begingroup\makeatletter\ifx\SetFigFont\undefined%
\gdef\SetFigFont#1#2#3#4#5{%
  \reset@font\fontsize{#1}{#2pt}%
  \fontfamily{#3}\fontseries{#4}\fontshape{#5}%
  \selectfont}%
\fi\endgroup%
\begin{picture}(2704,2912)(76,-2160)
\put( 76,-1909){\makebox(0,0)[lb]{\smash{\SetFigFont{12}{14.4}{\rmdefault}{\mddefault}{\updefault}$w^-$}}}
\put(1302,607){\makebox(0,0)[lb]{\smash{\SetFigFont{12}{14.4}{\rmdefault}{\mddefault}{\updefault}$w^+$}}}
\put(2780,-1626){\makebox(0,0)[lb]{\smash{\SetFigFont{12}{14.4}{\rmdefault}{\mddefault}{\updefault}$z^+$}}}
\put(1491,-1155){\makebox(0,0)[lb]{\smash{\SetFigFont{12}{14.4}{\rmdefault}{\mddefault}{\updefault}$2s$}}}
\put(1051,-1595){\makebox(0,0)[lb]{\smash{\SetFigFont{12}{14.4}{\rmdefault}{\mddefault}{\updefault}$z$}}}
\put(673,-966){\makebox(0,0)[lb]{\smash{\SetFigFont{12}{14.4}{\rmdefault}{\mddefault}{\updefault}$w$}}}
\put(1522,-305){\makebox(0,0)[lb]{\smash{\SetFigFont{12}{14.4}{\rmdefault}{\mddefault}{\updefault}$a_w$}}}
\put(2120,-1186){\makebox(0,0)[lb]{\smash{\SetFigFont{12}{14.4}{\rmdefault}{\mddefault}{\updefault}$a_z$}}}
\put(894,-588){\makebox(0,0)[lb]{\smash{\SetFigFont{12}{14.4}{\rmdefault}{\mddefault}{\updefault}$p$}}}
\put(1397,-1532){\makebox(0,0)[lb]{\smash{\SetFigFont{12}{14.4}{\rmdefault}{\mddefault}{\updefault}$q$}}}
\put(1806,-871){\makebox(0,0)[lb]{\smash{\SetFigFont{12}{14.4}{\rmdefault}{\mddefault}{\updefault}$r$}}}
\end{picture}